\documentclass[]{interact}

\usepackage{epstopdf}
\usepackage[caption=false]{subfig}
\usepackage[T1]{fontenc}
\usepackage{lmodern}
\usepackage{xcolor}
\usepackage{algorithm}
\usepackage{algorithmic}
\usepackage{url}
\usepackage{listings}
\usepackage{amsmath,amssymb,dsfont}

\usepackage[numbers,sort&compress]{natbib}
\bibpunct[, ]{[}{]}{,}{n}{,}{,}
\makeatletter
\def\NAT@def@citea{\def\@citea{\NAT@separator}}
\makeatother

\theoremstyle{plain}
\newtheorem{theorem}{Theorem}[section]
\newtheorem{lemma}[theorem]{Lemma}
\newtheorem{corollary}[theorem]{Corollary}
\newtheorem{proposition}[theorem]{Proposition}

\theoremstyle{definition}
\newtheorem{definition}[theorem]{Definition}
\newtheorem{example}[theorem]{Example}

\theoremstyle{remark}
\newtheorem{remark}{Remark}

\begin{document}

\articletype{}

\title{Structural Packing and Dyadic Factorization of Sparse Positive Definite Matrices}

\author{
\name{Micha{\l} Kos\textsuperscript{a}, Krzysztof Podg\'{o}rski\textsuperscript{b}, and Hanqing Wu\textsuperscript{b} \thanks{CONTACT Hanqing Wu Email: hanqing.wu@stat.lu.se}}
\affil{\textsuperscript{a}Institute of Mathematics, University of Wroc{\l}aw, Wroclaw, Poland; \\ \textsuperscript{b}Department of Statistics, Lund University, Lund, Sweden}
}

\maketitle

\begin{abstract}
    Efficient inversion of large sparse positive definite matrices requires exploiting sparsity patterns beyond those captured by conventional bandwidth reduction. In this work, we recast nested dissection, a prominent alternative, as a two-stage framework. The matrix was first packed into block-tridiagonal or dyadic form, followed by sparse Gram–Schmidt orthogonalization. This decomposition provided a unified perspective on sparse matrix factorization and inversion and identified dyadic structure as a fundamental component of sparse Cholesky factorization.
    For the first stage, we introduced a packing algorithm that recovered block-tridiagonal and dyadic patterns using a novel \(\ell_1\) criterion. Using approximate distances obtained through classical multidimensional scaling, the method was effective when the target structure was sufficiently represented among the nonzero entries. Iterative application could also remove structural noise and reveal hidden dyadic organization, corresponding to separator identification in nested dissection.
    For the second stage, we developed the theory of dyadically structured matrices. We derived sparse factorization and inversion procedures, analyzed their computational complexity, and obtained an efficient inversion algorithm. A modified version reduced the cost of inverting block-tridiagonal matrices, demonstrating the benefit of exploiting their structure directly rather than treating them as generic band matrices.
\end{abstract}


\section{Introduction}
\label{sec:intro}
Inverting large sparse matrices is a classical computational challenge in high-dimensional statistical analysis, spatial econometrics, and numerical solutions of algebraic and differential equations.
Matrix inversion is intrinsically related to classical algebraic techniques, such as least squares, Gaussian elimination, Gram–Schmidt orthogonalization, and Cholesky factorization.
We briefly account for these connections, which are often not spelled out.

Computing the inverse of a matrix can be reinterpreted as solving a linear system
$\mathbf y=\mathbf X \boldsymbol \beta$ for \(\boldsymbol\beta\).
If the solution is unique, then any method of finding it, such as Gaussian elimination, is equivalent to applying \(\mathbf X^{-1}\) to a given \(\mathbf y\). Our focus, computing the entire matrix \(\mathbf X^{-1}\), amounts to solving the system for the standard basis vectors.
The method of least squares extends this viewpoint to settings where an exact solution is unavailable, replacing it with $\widehat{\boldsymbol \beta}=
    \left(\mathbf X^\top \mathbf X\right)^{-1}\mathbf X^\top \mathbf y$ which minimizes the sum of squares of the entries of $\mathbf y-\mathbf X\boldsymbol \beta$.
In this case, the problem is transferred to inverting the positive definite Gram matrix $\boldsymbol \Sigma = \mathbf X^\top \mathbf X$. If $\mathbf X$ is sparse, then \(\boldsymbol\Sigma\) may, but need not, be sparse.


Since $\widehat{\boldsymbol \beta}$ is obtained by minimizing the squared Euclidean distance, $\mathbf X \widehat{\boldsymbol \beta}$ represents the projection of $\mathbf y$ onto the column space of $\mathbf X$. Such a projection can be obtained by orthonormalizing the columns of $\mathbf{X}$, which can be formally stated as finding a matrix $\mathbf P$ such that
$
    \mathbf X\mathbf P
$
has orthonormal columns, i.e.,
\begin{equation}
    \label{eq:ON}
    \mathbf P^\top\boldsymbol \Sigma \mathbf P=\mathbf I.
\end{equation}
The direct consequences of \eqref{eq:ON} are $\mathbf P^{-1} =\mathbf P^\top \boldsymbol \Sigma$,  $\boldsymbol \Sigma^{-1} =\mathbf P \mathbf P^\top$, and thus
\begin{align}
    \label{eq:PreCholesky}
    \boldsymbol \Sigma & ={\mathbf P^{-1}}^\top \mathbf P^{-1}, \\
    \label{eq:pre-QR}
    \mathbf X          & =\mathbf Q \mathbf P^{-1},
\end{align}
where $\mathbf Q=\mathbf X \mathbf P$ has orthonormal columns. Consequently, once $\mathbf P$ is found, the inversion, projection, factorization, and least-squares solution are obtained. We emphasize that, although $\mathbf P$ is generally not unique, its chosen form is crucial in the computational aspect of these operations.
When $\mathbf{P}$ is computed via the Gram-Schmidt process, both  $\mathbf{P}$ and its inverse are triangular. In this case, \eqref{eq:PreCholesky} corresponds to Cholesky factorization and \eqref{eq:pre-QR} to QR decomposition; indeed, Gram-Schmidt orthogonalization can be interpreted as Gaussian elimination applied to $\boldsymbol \Sigma$, see \cite{PursellT1991}.

Choosing the structure of $\mathbf P$ properly when $\boldsymbol{\Sigma}$ is sparse is therefore computationally important. A prominent method is nested dissection \cite{Alan1973}. There is a substantial body of work in which nested dissection is studied in the context of specific graph classes, such as planar graphs, finite element meshes, and graph families admitting balanced separators \cite{Alan1973,Alan1978,LiptonRT1979,Gilbert:1986aa,CambierChB2020}, often through algorithms tailored to particular applications.
The objective of the present work is complementary. Rather than focusing on highly optimized, problem-specific implementations, we seek to identify the structural principles underlying sparse factorization and inversion. This leads to a unified interpretation of nested dissection and to algorithms that emphasize conceptual clarity and generality. Extensions tailored to specific settings can be naturally developed within our proposed framework, and we leave them for future work.

In this work, we interpret sparse matrix factorization and inversion based on nested dissection as a two-step procedure.
The first step recursively finds the {\it separators} and then constructs a permutation matrix $\boldsymbol \pi$ such that $\boldsymbol \pi^\top \boldsymbol \Sigma \boldsymbol \pi$ has a {\it dyadic form}.
The second step applies dyadic Gram-Schmidt orthogonalization to the reordered matrix obtained from the first, resulting in a sparse factorization. In George's original formulation \cite{Alan1973}, the separator hierarchy determines the elimination ordering used to obtain an efficient sparse triangular decomposition.
The key benefit of nested dissection is that the ordering limits {\em fill-in}, that is, nonzero entries appearing in the factor matrices at positions that are zero in the reordered matrix.

Our main contributions arise from a dyadic formulation of nested dissection. We interpret matrices amenable to nested dissection as having a dyadic sparsity pattern up to permutation. Within this class, two main challenges arise: first, the construction of an appropriate dyadic reordering (packing); and second, the sparse factorization and inversion of the resulting dyadic matrix. To address the first challenge, we introduce a generic structural packing method for sparse matrices. It employs multidimensional scaling and aims to minimize the $\ell_1$-norm of the half-width vector. The method is not tailored to a specific graph class and therefore applies more broadly. Compared with the classical bandwidth minimization based on the \(\ell_\infty\)-norm, our algorithm based on the \(\ell_1\) criterion is less sensitive to a small number of nonzero entries far from the diagonal. Moreover, the packing algorithm can be applied recursively to unveil hidden dyadic structure, corresponding to successive separator identification. We tackle the second challenge by an efficient dyadic Gram–Schmidt orthogonalization, for which we explicitly evaluate the computational cost. This procedure yields a dyadic factorization of the resulting matrix. We derive computational complexity results for both dyadic factorization and inversion, showing in particular that the factorization has quasi-linear complexity.

For expository clarity, we present the two contributions mentioned above in reverse order. The paper is organized as follows. In Section~\ref{sec:dyadmat}, we discuss the algebraic properties of dyadic matrices.
Section~\ref{sec:dyadalg} contains the computational complexity results for the dyadic factorization. The novel packing algorithm is presented in Section~\ref{sec:packing}, and the corresponding simulation results are given in Section~\ref{sec:sim}. We provide additional technical details and auxiliary information in the Appendix. Finally, an R implementation of the dyadic matrix algebra, including the dyadic factorization algorithm, is available in the R package \texttt{DyadiCarma} \cite{KosPW2025}; an implementation of the packing algorithm can be found in the GitHub repository \cite{packing_algorithm_github} together with code reproducing all simulation and computational results of the paper.

\section{Dyadic matrices, their decompositions and inversions}
\label{sec:dyadmat}
We discuss several subalgebras of the associative algebra $\mathcal M$ of $d\times d$ matrices over $\mathbb R$.
These subalgebras are characterized by dyadic structures as described below.

Let $k, s, r  \in \mathbb{N}$ and let $\mathcal{A} \subset [s] \times [r]$ be a set of indices.
Here and in what follows $[s]=\{1,\dots, s\}$.
Let a $ks \times kr$ matrix $\mathbf{A}$ be partitioned into a $s \times r$ matrix of $k \times k$ blocks
$$
    \mathbf{A} = \begin{bmatrix}
        \mathbf{A}_{1,1} & \cdots & \mathbf{A}_{1,r} \\
        \vdots           & \ddots & \vdots           \\
        \mathbf{A}_{s,1} & \cdots & \mathbf{A}_{s,r}
    \end{bmatrix}.
$$
We say that $\mathbf A$ has a (regular) block-sparsity pattern $
    \mathcal{A},
$  if $\mathbf{A}_{i,j} = \mathbf{0}_{k \times k}$ for $(i,j) \notin \mathcal{A}$, where $\mathbf 0_{k \times k}$ stands for the $k\times k$ zero matrix. Throughout this work, the block size $k$ is fixed and omitted from the notation.
We abbreviate this as $\mathbf{A} \in \mathcal{A}$.
\begin{remark}
    The concept of a block-sparsity pattern naturally extends to irregular cases, where matrices are partitioned into blocks of varying sizes or where the blocks are not necessarily square. However, for simplicity, we restrict our discussion to the regular case where all blocks are square and have the same size $k$. Although we do not consider this more general setup here, any irregular or incomplete dyadic matrix can be naturally embedded into a regular and complete one by filling missing diagonal blocks with identity matrices and the non-diagonal ones with zero matrices.
    The algebraic results and corresponding algorithms extend through this embedding, with their complexity expressed in terms of the augmented dimension and block size.
    Such generalizations would be relevant for practical implementation and for connecting the results to nested dissection algorithms.
    Their mathematical derivation would follow the arguments presented here, with appropriate adjustments to the notation and statements.
\end{remark}

In what follows, for a partitioned $ks\times kr$ matrix $\mathbf{A}$ with the corresponding sparsity pattern $\mathcal{A}$, we consider the mapping $
    \mathbf{A}_{\mathcal{A}}:(i,j)\mapsto \mathbf{A}_{i,j}, (i,j)\in\mathcal{A},
$ that assigns to each index pair its corresponding $k\times k$ block.
Similarly, if $\mathbf{A}$ is an $m\times kr$ matrix partitioned into $r$ consecutive blocks of $k$ columns with $\mathcal{B}\subset[r]$, we define
$
    \mathbf{A}_{\cdot,\mathcal{B}}: j\mapsto \mathbf{A}_{\cdot,j}, j\in\mathcal{B},
$
as the mapping that assigns to each column block index its corresponding $m\times k$ submatrix. Similarly, we define the row-wise analog $ \mathbf{A}_{\mathcal{B},\cdot}$.

As noted, dyadically structured matrices can be defined for an arbitrary $d$. In this work, we focus exclusively on the regular completely dyadic $d\times d$ matrices for which  $d=k\left(2^N-1\right)$, where $N,k\in \mathbb N$, with $k$ fixed throughout our discussion.
We refer to $N$ as the {\em height} of a dyadic structure.
Let us define the pyramids $\mathcal P(N)$, $\mathcal P'(N)$ of height $N$, consisting of sets of integers referred to as levels
\begin{align}
    \label{eq:dyadord}
    I_{l}  & = \bigcup_{r=1}^{2^{N-l}} I_{r, l},\,\, I_{r, l} = \left\{2^{l-1}(2r-1)\right\}, \ l=1,\dots,N, \\
    \label{eq:dyadInt}
    I'_{l} & = \bigcup_{r=1}^{2^{N-l}} I'_{r, l},\,\, I'_{r, 1}  = I_{r, 1},\,\,
    I'_{r, l}= I_{r, l} \cup I_{2 r-1, l-1}^{\prime} \cup I_{2 r, l-1}^{\prime},\, l=2, \ldots, N.
\end{align}

There is an explicit decomposition of $I'_{r,l}$ in terms of $I_{r,l}$. For $l\in \{1,\dots, N\}$, the index set $I'_{r,l}$, $r\in \{1,\dots, 2^{N-l}\}$ is made of consecutive integers
$2^{l-1}(2 r - 1) + \left\{\pm 2^{l-1} + 1, \dots, \pm 1, 0 \right\}.
$
This decomposition, which is a direct consequence of \eqref{eq:dyadInt}, is useful for efficient implementation of the algebra of dyadic matrices.

We define the {\em horizontally dyadic structure} as the block-sparsity pattern given by
$$
    \mathcal{HD}(N) = \bigcup_{l=1}^N\bigcup_{r=1}^{2^{N-l}} \mathbf I_{r, l} \subset [2^N-1]^2,
$$
where
$
    \mathbf I_{r, l} = I_{r, l} \times I'_{r, l}.
$
Similarly, let
$
    \mathbf{I}_{r, l}^{\top}=\left\{(i, j) \in[2^N-1]^2:(j, i) \in \mathbf{I}_{r, l}\right\}.
$
The {\em vertically dyadic structure}, denoted by $\mathcal{V} \mathcal{D}(N)$, is defined as the block-sparsity pattern given by the union of $\mathbf I_{r,l}^\top$'s.
Finally, the {\em symmetrically dyadic structure} $\mathcal{SD}(N)$ is defined as the union of
$
    \mathbf{J}_{r, l}=\mathbf{I}_{r, l} \cup \mathbf{I}_{r, l}^{\top}
$.

A $d\times d$ matrix is called horizontally, vertically, or symmetrically dyadic if it has a block-sparsity pattern corresponding to $\mathcal{HD}(N)$, $\mathcal{VD}(N)$, or $\mathcal{SD}(N)$, respectively. For brevity, the corresponding classes of matrices are also denoted by $\mathcal{HD}(N)$, $\mathcal{VD}(N)$, and $\mathcal{SD}(N)$.
The dyadic structures are identified with the respective matrices $\mathbf{H}(N)$, $\mathbf{V}(N)$, and $\mathbf{S}(N)$, which have $1$’s in all positions corresponding to potentially nonzero blocks. When the value of \(k\) matters, we use $\mathbf{H}(N, k)$, $\mathbf{V}(N, k)$, and $\mathbf{S}(N, k)$ to represent the corresponding 0-1 matrix with block size \(k\).

For $N=4$, the three block-sparsity patterns are shown through $*$ in the following matrices.
    {\tiny
        \setlength\arraycolsep{2pt}
        \begin{equation*}
            \begin{bmatrix}
                *               &                 &                 &                 &                 &                 &                 &                &                 &                 &                 &                 &                 &                 &                 \\
                \color{red} *   & \color{red} *   & \color{red} *   &                 &                 &                 &                 &                &                 &                 &                 &                 &                 &                 &                 \\
                                &                 & *               &                 &                 &                 &                 &                &                 &                 &                 &                 &                 &                 &                 \\
                \color{green} * & \color{green} * & \color{green} * & \color{green} * & \color{green} * & \color{green} * & \color{green} * &                &                 &                 &                 &                 &                 &                 &                 \\
                                &                 &                 &                 & *               &                 &                 &                &                 &                 &                 &                 &                 &                 &                 \\
                                &                 &                 &                 & \color{red} *   & \color{red} *   & \color{red} *   &                &                 &                 &                 &                 &                 &                 &                 \\
                                &                 &                 &                 &                 &                 & *               &                &                 &                 &                 &                 &                 &                 &                 \\
                \color{blue} *  & \color{blue} *  & \color{blue} *  & \color{blue} *  & \color{blue} *  & \color{blue} *  & \color{blue} *  & \color{blue} * & \color{blue} *  & \color{blue} *  & \color{blue} *  & \color{blue} *  & \color{blue} *  & \color{blue} *  & \color{blue} *  \\
                                &                 &                 &                 &                 &                 &                 &                & *               &                 &                 &                 &                 &                 &                 \\
                                &                 &                 &                 &                 &                 &                 &                & \color{red} *   & \color{red} *   & \color{red} *   &                 &                 &                 &                 \\
                                &                 &                 &                 &                 &                 &                 &                &                 &                 & *               &                 &                 &                 &                 \\
                                &                 &                 &                 &                 &                 &                 &                & \color{green} * & \color{green} * & \color{green} * & \color{green} * & \color{green} * & \color{green} * & \color{green} * \\
                                &                 &                 &                 &                 &                 &                 &                &                 &                 &                 &                 & *               &                 &                 \\
                                &                 &                 &                 &                 &                 &                 &                &                 &                 &                 &                 & \color{red} *   & \color{red} *   & \color{red} *   \\
                                &                 &                 &                 &                 &                 &                 &                &                 &                 &                 &                 &                 &                 & *               \\
            \end{bmatrix}
            \begin{bmatrix}
                * & \color{red} * &   & \color{green} * &   &               &   & \color{blue} * &   &               &   &                 &   &               &   \\
                  & \color{red} * &   & \color{green} * &   &               &   & \color{blue} * &   &               &   &                 &   &               &   \\
                  & \color{red} * & * & \color{green} * &   &               &   & \color{blue} * &   &               &   &                 &   &               &   \\
                  &               &   & \color{green} * &   &               &   & \color{blue} * &   &               &   &                 &   &               &   \\
                  &               &   & \color{green} * & * & \color{red} * &   & \color{blue} * &   &               &   &                 &   &               &   \\
                  &               &   & \color{green} * &   & \color{red} * &   & \color{blue} * &   &               &   &                 &   &               &   \\
                  &               &   & \color{green} * &   & \color{red} * & * & \color{blue} * &   &               &   &                 &   &               &   \\
                  &               &   &                 &   &               &   & \color{blue} * &   &               &   &                 &   &               &   \\
                  &               &   &                 &   &               &   & \color{blue} * & * & \color{red} * &   & \color{green} * &   &               &   \\
                  &               &   &                 &   &               &   & \color{blue} * &   & \color{red} * &   & \color{green} * &   &               &   \\
                  &               &   &                 &   &               &   & \color{blue} * &   & \color{red} * & * & \color{green} * &   &               &   \\
                  &               &   &                 &   &               &   & \color{blue} * &   &               &   & \color{green} * &   &               &   \\
                  &               &   &                 &   &               &   & \color{blue} * &   &               &   & \color{green} * & * & \color{red} * &   \\
                  &               &   &                 &   &               &   & \color{blue} * &   &               &   & \color{green} * &   & \color{red} * &   \\
                  &               &   &                 &   &               &   & \color{blue} * &   &               &   & \color{green} * &   & \color{red} * & * \\
            \end{bmatrix}
            \begin{bmatrix}
                *               & {\color{red} *}   &                 & {\color{green}
                *}              &                   &                 &                 & {\color{blue}
                *}              &                   &                 &                 &                   &                 &                 &                 \\
                {\color{red} *} & {\color{red} *}   & {\color{red} *} & {\color{green}
                *}              &                   &                 &                 & {\color{blue}
                *}              &                   &                 &                 &                   &                 &                 &                 \\
                                & {\color{red} *}   & *               & {\color{green}
                *}              &                   &                 &                 & {\color{blue}
                *}              &                   &                 &                 &                   &                 &                 &                 \\
                {\color{green}
                *}              & {\color{green}
                *}              & {\color{green}
                *}              & {\color{green}
                *}              & {\color{green}
                *}              & {\color{green}
                *}              & {\color{green}
                *}              & {\color{blue}
                *}              &                   &                 &                 &                   &                 &                 &                 \\
                                &                   &                 & {\color{green}
                *}              & *                 & {\color{red} *} &                 & {\color{blue}
                *}              &                   &                 &                 &                   &                 &                 &                 \\
                                &                   &                 & {\color{green}
                *}              & {\color{red} *}   & {\color{red} *} & {\color{red} *} & {\color{blue}
                *}              &                   &                 &                 &                   &                 &                 &                 \\
                                &                   &                 & {\color{green}
                *}              &                   & {\color{red} *} & *               & {\color{blue}
                *}              &                   &                 &                 &                   &                 &                 &                 \\ {\color{blue}
                *}              & {\color{blue}
                *}              & {\color{blue}
                *}              & {\color{blue}
                *}              & {\color{blue}
                *}              & {\color{blue}
                *}              & {\color{blue}
                *}              & {\color{blue}
                *}              & {\color{blue}
                *}              & {\color{blue}
                *}              & {\color{blue}
                *}              & {\color{blue}
                *}              & {\color{blue}
                *}              & {\color{blue}
                *}              & {\color{blue}
                *}                                                                                                                                                \\
                                &                   &                 &                 &                   &                 &                 & {\color{blue}
                *}              & *                 & {\color{red} *} &                 & {\color{green} *} &                 &                 &                 \\
                                &                   &                 &                 &                   &                 &                 & {\color{blue}
                *}              & {\color{red} *}   & {\color{red} *} & {\color{red} *} & {\color{green} *} &                 &                 &                 \\
                                &                   &                 &                 &                   &                 &                 & {\color{blue}
                *}              &                   & {\color{red} *} & *               & {\color{green} *} &                 &                 &                 \\
                                &                   &                 &                 &                   &                 &                 & {\color{blue}
                *}              & {\color{green}
                *}              & {\color{green}
                *}              & {\color{green}
                *}              & {\color{green}
                *}              & {\color{green}
                *}              & {\color{green}
                *}              & {\color{green} *}                                                                                                               \\
                                &                   &                 &                 &                   &                 &                 & {\color{blue}
                *}              &                   &                 &                 & {\color{green} *} & *               & {\color{red} *} &                 \\
                                &                   &                 &                 &                   &                 &                 & {\color{blue}
                *}              &                   &                 &                 & {\color{green} *} & {\color{red} *} & {\color{red} *} & {\color{red} *} \\
                                &                   &                 &                 &                   &                 &                 & {\color{blue}
                *}              &                   &                 &                 & {\color{green} *} &                 & {\color{red} *} & *               \\
            \end{bmatrix}
        \end{equation*}
    }
To see the different dyadic structure levels, they are marked with different colors: the 1st level - {\it black}, the 2nd level - {\it red}, the 3rd level - {\it green}, and the fourth level - {\it blue}. For $\mathbf S=\mathbf S(N,k)$, $\mathbf S_{\mathbf I_{1,1}} = \mathbf 1_{k\times k}$ and $\mathbf S_{\mathbf J_{1,2}}$ is made of five $k\times k$ matrices forming a cross marked in {\it red} at the top-left corner of the right-hand side matrix.

\begin{theorem}
    \label{th:alg} The classes $\mathcal{HD}(N)$, $\mathcal{VD}(N)$, $\mathcal{SD}(N)$ are linear spaces, and the first two are (associative) subalgebras of $\mathcal M$. Moreover,
    \begin{align*}
        \mathcal{HD}(N)\mathcal{VD}(N)
                                       & \subseteq\mathcal{SD}(N)=\mathcal{HD}(N)+\mathcal{VD}(N), \\
        \mathcal{HD}(N)\mathcal{SD}(N) & =
        \mathcal{SD}(N)\mathcal{VD}(N)= \mathcal{SD}(N).
    \end{align*}
\end{theorem}
The assertions follow directly from the nesting relations among the index sets \(I'_{r,l}\), and the proof is omitted.
\begin{remark}
    It should be noted that $\mathcal{VD}(N)\mathcal{HD}(N)$ need not be dyadic in any sense.
\end{remark}
The next result emphasizes the recurrent nature of the dyadic structures and connects them to the Gramian representation of PD dyadic matrices.
\begin{proposition}
    \label{prop:propdyad} The following properties hold
    \begin{description}
        \item[\it i)]
              Let $\boldsymbol \Sigma=\mathbf F^\top \mathbf F\in \mathcal{SD}(N)$ and define $\mathbf F_{r,l}=\mathbf F_{\cdot, I_{r, l}}$, $r=1,\dots, 2^{N-l}$, $l=1,\dots N$.
              If
              $I'_{r,l}\cap I'_{r',l'}=\emptyset$, then the linear spaces spanned by columns of $\mathbf F_{r,l}$ and $\mathbf F_{r',l'}$ are orthogonal.

        \item[\it ii)] We have the following recursive relationship

                  {
                      \begin{equation*}
                          \begin{aligned}
                              \mathbf H(N+1) & = \begin{bmatrix}
                                                     \mathbf H(N)           & \mathbf 0_{d\times k} & \mathbf 0_{d\times d} \\
                                                     \mathbf 1_{d \times d} & \mathbf 1_{k\times k} & \mathbf 1_{k\times d} \\
                                                     \mathbf 0_{d\times d}  & \mathbf 0_{d\times k} & \mathbf H(N)
                                                 \end{bmatrix},\;
                              \mathbf V(N+1)   = \begin{bmatrix}
                                                     \mathbf V(N)           & \mathbf 1_{d\times k} & \mathbf 0_{d\times d} \\
                                                     \mathbf 0_{k \times d} & \mathbf 1_{k\times k} & \mathbf 0_{k\times d} \\
                                                     \mathbf 0_{d\times d}  & \mathbf 1_{d\times k} & \mathbf V(N)
                                                 \end{bmatrix}, \; \\
                              \mathbf S(N+1) & = \begin{bmatrix}
                                                     \mathbf S(N)           & \mathbf 1_{d\times k} & \mathbf 0_{d\times d} \\
                                                     \mathbf 1_{k \times d} & \mathbf 1_{k\times k} & \mathbf 1_{k\times d} \\
                                                     \mathbf 0_{d\times d}  & \mathbf 1_{d\times k} & \mathbf S(N)
                                                 \end{bmatrix}.
                          \end{aligned}
                      \end{equation*}
                  }

        \item[\it iii)] If we remove the rows and columns from a dyadic matrix of height \( N \) that correspond to the index sets \( I_{r,1} \) for \( r=1,\dots,2^{N-1} \), the resulting matrix remains dyadic but with a reduced height of \( N-1 \).
    \end{description}
\end{proposition}

If $\mathbf A$ is a dyadic matrix of any type, then we identify it with the family of submatrices specified by the index set of dyadic structures $\mathbf A_{\mathcal I}$.
Then we do not distinguish between $\mathbf A$ and $\mathbf A_{\mathcal I}$.
The operation given in {\it iii)} is defined as a subsampling along $\check{\mathcal I}= \mathcal{HD}(N-1)$ and the resulting symmetric dyadic matrix is written as $\mathbf A_{\check{\mathcal I}}$.
This extends to any $\check{\mathcal I}=\mathcal{HD}(M)$ for any $M<N$.

In the next theorem, we present factorization results for dyadic PD matrices.
Since PD matrices are symmetric, only the symmetrically dyadic class is relevant and we hence omit the reference to symmetry.
\begin{theorem}
    \label{th:factor}
    Let $\mathbf \Sigma \in \mathcal{SD}(N)$ be positive definite. Then there exists $\mathbf P
        \in \mathcal{VD}(N)$ such that
    $
        \mathbf P^\top \boldsymbol{\Sigma}\mathbf P=\mathbf I
    $.
    Moreover, $\mathbf R=\mathbf P^{-1}=\mathbf P^\top \boldsymbol \Sigma \in \mathcal{VD}(N)$.
\end{theorem}
\begin{corollary}
    With the notation of Theorem~\ref{th:factor},
    $\boldsymbol \Sigma=\mathbf R^\top\mathbf R,$ $
        \boldsymbol \Sigma^{-1}  = \mathbf P \mathbf P^\top$.
\end{corollary}
The algorithm delivering $ \mathbf {P} $ is presented in the next section, and its justification serves as a proof of the theorem, see Lemma~\ref{lem:IHDED} and Theorem~\ref{th:maindyadic}.

\begin{remark}
    It should be noted that $\boldsymbol \Sigma^{-1}$ is typically not a dyadic matrix.
\end{remark}

\section{Dyadic algorithm and its complexity}
\label{sec:dyadalg}
For a PD $\boldsymbol{\Sigma} \in \mathcal{SD}(N)$, the algorithm returns $\mathbf{P} \in \mathcal{VD}(N)$ satisfying
$
    \mathbf{P}^{\top} \boldsymbol{\Sigma} \mathbf{P}=\mathbf{I}
$.
Since the output is vertically dyadic, the algorithm computes the blocks $\mathbf{P}_{\mathcal{I}}$ with $\mathcal{I}=\mathcal{VD}(N)$, which uniquely determine $\mathbf{P}$.
The algorithm proceeds through levels $l=1, \ldots, N$, computing these blocks through an iterative procedure. Because this method specifically utilizes the dyadic structures, we refer to it as the dyadic algorithm, although its overall scheme is a general sequential orthogonalization, algebraically equivalent to sparse Gaussian elimination or sparse Cholesky decomposition.

In the process, an unspecified orthogonalization $\mathbf G(\overline{\boldsymbol \Sigma})$ of a $k\times k$ PD matrix $\overline{\boldsymbol \Sigma}$ is utilized, and for this one can consider any generic algorithm that orthogonalizes a matrix $\overline{\mathbf F}$ of $k$ columns such that $
    \overline{\boldsymbol \Sigma}=
    \overline{\mathbf F}^\top \overline{\mathbf F}$.
The matrix $\mathbf P$ provided by the algorithm is uniquely determined once $\mathbf G(\cdot)$ is specified.
A simple choice for $\mathbf G(\cdot)$ is the classical Gram-Schmidt method. In general, the entire idea for orthogonalizing a dyadic matrix is similar to the Gram-Schmidt method; thus, the process of obtaining $\mathbf P$ can also be naturally referred to as the dyadically structured Gram-Schmidt method.

Implementations of the dyadic algorithm for general symmetric dyadic matrices, an accelerated version for block-tridiagonal matrices, and basic algebraic operations for dyadic matrices are provided in the R package \texttt{DyadiCarma} \cite{KosPW2025}. The package leverages {\tt Rcpp} and {\tt RcppArmadillo} to deliver high-performance computation.

To evaluate the complexity of the algorithm, we count the number of floating-point operations (flops). Typically, an addition, a multiplication, or a comparison is counted as one flop, while divisions, square roots, and exponentials are assigned a higher but constant cost \cite{hockney_science_1996}. However, following \cite{trefethen_numerical_1997}, we treat each of these operations as a single flop since this simplification does not affect the overall complexity in the big-$O$ notation. A similar complexity result was obtained in earlier works in the context of sparse Gaussian elimination \cite{rose1972graph,LiptonRT1979}, providing corroboration for the optimality of the dyadic algorithm.

The dyadic algorithm uses sequential orthogonalization, which is an extension of the Gram-Schmidt method.
In fact, sequential orthogonalization, described in Algorithm~\ref{alg:dyadic}, works for any $\mathbf F=\begin{bmatrix}\mathbf F_1&\dots &\mathbf F_N\end{bmatrix}$ and leads to
$\mathbf P$ in \eqref{eq:ON} for $\boldsymbol \Sigma=\mathbf F^\top \mathbf F$.
\begin{algorithm}
    \caption{Sequential orthogonalization $\mathbf P(\boldsymbol \Sigma)$}
    \label{alg:dyadic}
    \begin{algorithmic}\vspace{2mm}

        \STATE{{\sc Input:} $\mathbf \Sigma$ corresponding to a $d\times d$ symmetrically dyadic matrix}\vspace{1mm}
        \STATE{$\mathbf P=\mathbf G_1(\boldsymbol \Sigma_{1,1})$};  \COMMENT{{\tt orthonormalization of $\mathbf F_1$, i.e., $\mathbf P^\top \mathbf F_1^\top \mathbf F_1 \mathbf P = \mathbf I$
                }}\vspace{1mm}
        \FOR{$l=2$ to $N$}\vspace{2mm}
        \STATE{$\check{\mathbf \Sigma} = \begin{bmatrix}
                    \boldsymbol \Sigma_{1,l} \\
                    \vdots                   \\
                    \boldsymbol \Sigma_{l-1,l}
                \end{bmatrix}$};\vspace{1mm}
        \STATE{$\check{\boldsymbol \Sigma}^\prime = {\mathbf P}^\top\check{\boldsymbol\Sigma}$;\, $\tilde{\boldsymbol \Sigma} = \boldsymbol \Sigma_{l,l} - \check{\boldsymbol \Sigma}^\prime{}^\top \check{\boldsymbol \Sigma}^\prime$}; \vspace{1mm}
        \STATE{$\mathbf A=\mathbf P \check{\boldsymbol \Sigma}^\prime$;\,${\mathbf G} = \mathbf G_l(\tilde{\boldsymbol \Sigma})$};\vspace{1mm}
        \STATE{${\mathbf P} = \begin{bmatrix}\mathbf P & -\mathbf A \mathbf G \\
               \mathbf 0 & \mathbf G\end{bmatrix}$}; \COMMENT{{\tt  Orthonormalization of $\begin{bmatrix}\mathbf F_1&\dots &\mathbf F_l\end{bmatrix}$}}\vspace{1mm}
        \ENDFOR \vspace{1mm}
        \STATE{{\sc Return:} $\mathbf P$ orthogonalizing $\boldsymbol \Sigma$, as in factorization corresponding to a vertically dyadic matrix.}
    \end{algorithmic}
\end{algorithm}
We assume that each Gramian $\boldsymbol \Sigma_{l,l}=\mathbf F_l^\top \mathbf F_l$ possesses a structure that allows for a generic orthonormalization procedure denoted by $\mathbf G_l(\cdot)$, where the argument is a PD matrix that has a desired structure.
The algorithm produces a matrix $\mathbf P$ that satisfies \eqref{eq:ON}, as long as the matrices $\tilde{\boldsymbol \Sigma}$ in the loop retain the structure required for the application of each $\mathbf G_l$, and its complexity is summarized in the following result.

\begin{proposition}
    \label{prop:complexity}
    Let $g_l$ be the complexity of evaluating $\mathbf G_l(\tilde{ \boldsymbol \Sigma})$ for a specific matrix $\tilde{\boldsymbol \Sigma}$ and $a_l$, $b_l$, $c_l$, $d_l$ are the costs of evaluating ${\mathbf P}^\top\check{\mathbf \Sigma}$, $\mathbf P \check{\boldsymbol \Sigma}^\prime$, $\boldsymbol \Sigma_{l,l} - \check{\mathbf \Sigma}^\top{}^\prime \check{\mathbf \Sigma}^\prime$, $\mathbf A \mathbf G$, respectively, in the $l$th step of the algorithm, then the total cost of the sequential orthogonalization algorithm is
    $
        \sum_{l=1}^Ng_l+\sum_{l=2}^N\left(a_l+b_l+c_l+d_l\right)
    $.
\end{proposition}

The main algorithm of this work utilizes sequential orthogonalization in the context of symmetric dyadic matrices $\boldsymbol \Sigma$ that are the Gramians of
$
    \mathbf F=\begin{bmatrix}\mathbf F_1&\dots &\mathbf F_N\end{bmatrix}$,
where
$\mathbf F_l=\begin{bmatrix}\mathbf F_{1,l}&\dots & \mathbf F_{2^{N-l},l}\end{bmatrix}$, $l=1,\dots,N$ but considered in the dyadic order.
The transfer from the sequential order to the dyadic one is represented by
$$
    \{ 1,\dots, 2^N-1\} \ni i=2^N(1-2^{-l+1})+r \longleftrightarrow (r,l) \in \bigcup_{l=1}^N\{1,\dots, 2^{N-l}\}\times \{l\}.
$$
Thus, the sequential orthogonalization for the dyadic matrices $\boldsymbol \Sigma\in \mathcal{SD}(N)$ becomes dyadic one by utilizing one-to-one mapping from vectors in the sequential order to the dyadic order expressed in levels $l\in \{1,\dots, N\}$ and locations $r\in \{1,\dots,2^{N-l}\}$ within the level $l$.
This dyadic pyramid can be turned into a linear one using \eqref{eq:dyadord}.

The well-established numerical complexity of the Gram-Schmidt orthogonalization method applied to a $k\times k$ Gramian matrix is expressed by $O(k^3)$ flops.
In the sequential algorithm that runs through levels of the dyadic pyramid, orthogonalization $\boldsymbol G_l(\cdot)$ is based on Gram-Schmidt orthogonalizations of $2^{N-l}$ $k\times k$ submatrices from which $\tilde{\boldsymbol \Sigma}$ is formed at the diagonal (and zero otherwise).
We conclude that in the complexity of Proposition~\ref{prop:complexity}:
$$
    \sum_{l=1}^N g_l=\sum_{l=1}^N 2^{N-l}O(k^3)=(2^N-1)O(k^3) = O(dk^2).
$$
Thus, it suffices to discuss the remaining terms in the complexity formula of Proposition~\ref{prop:complexity}. With the same transfer scheme, matrices $\check{\mathbf \Sigma}$ and $\mathbf P$ during the sequential orthogonalization have specific block-sparsity patterns, namely $\mathcal{ED}(N, l)$ and $\mathcal{IVD}(N, l)$.

We consider the {\it diagonal block-sparsity pattern} of block size $k$ for matrices of size $k2^{N-l}\times k2^{N-l}$, $l = 1, \dots, N$:
\begin{equation}
    \mathcal{D}(N, l) = \{(i, j):  i \in \left[ 2^{N-l}\right], j = i \},
\end{equation}
where, as before, $[s] = \{1,\dots,s\}$.
Furthermore, the {\it elongated dyadic diagonal block-sparsity pattern} of block size $k$ is defined for matrices of size $k\left(2^{N} - 2^{N-l+1}\right) \times k2^{N-l}$ , $l=2,\dots, N$, by
\begin{equation}
    \label{eq:ED}
    \mathcal {ED}(N,l)=
    \left\{
    (i, j):  j \in [2^{N-l}],
    i \in
    (j-1)\left(2^{l}-2\right)+[2^{l}-2]
    \right \}.
\end{equation}
This pattern can be interpreted as the block-sparsity pattern of the submatrix
\begin{equation}
    \label{eq:edd_full}
    \mathbf \Sigma_{\mathcal{I}_{\mathcal{ED}}}, \text{ where } \mathcal{I}_{\mathcal{ED}} = \bigcup_{r=1}^{2^{N-l}} \left(\mathbf I_{r,l}^\top\setminus (I_{r,l}\times I_{r,l})\right) \text{ and }\mathbf \Sigma \in \mathcal{SD}(N)
\end{equation}
Define
$
    \mathcal {ED}(N, l)_{r} =  \{(i, r):  i \in  (r - 1) (2^l - 2) +  [2^l - 2]\}$,
which identifies the elongated diagonal blocks indexed by $\mathbf I_{r,l}^\top\setminus (I_{r,l}\times I_{r,l})$ in the corresponding  $\mathbf \Sigma$ in \eqref{eq:edd_full}. If matrix $\mathbf E\in \mathcal{ED}(N, l)$, we can decompose it as
\begin{equation}
    \label{eq:decompE}
    \mathbf E = \sum_{r=1}^{2^{N-l}} \mathbf E_r,
\end{equation}
where $\mathbf E_r$ has the block-sparsity pattern $\mathcal {ED}(N, l)_{r}$.

\begin{example}
    With the notation of Proposition~\ref{prop:propdyad}, if a PD dyadic matrix $\boldsymbol \Sigma=\mathbf F^\top \mathbf F$, then the submatrix
    $
        \mathbf F_1^\top \mathbf F_2
    $
    is $\mathcal {ED}(N,l)$-elongated diagonal, where
    \begin{align*}
        \mathbf F_1=\begin{bmatrix}
                        \mathbf F_{r,l'}
                    \end{bmatrix}_{r=1,l'=1}^{2^{N-l},l-1},\,\,\,
        \mathbf F_2 & =\begin{bmatrix}
                           \mathbf F_{r,l}
                       \end{bmatrix}_{r=1}^{2^{N-l}}.
    \end{align*}
    The block matrices ${\begin{bmatrix}
            \mathbf F_{r,l'}
        \end{bmatrix}_{l'=1}^{l-1}}^\top\mathbf F_{r,l}$ are the elongated diagonal blocks corresponding to $\mathcal {ED}(N, l)_r$, $r=1,\dots, 2^{N-l}$.
\end{example}

Similarly, the incomplete dyadic block-sparsity pattern can be conveniently characterized by submatrices of $\mathbf \Sigma$ as in~\eqref{eq:edd_full}.
For matrices of size $k\left(2^N-2^{N-l+1}\right)\times k\left(2^N-2^{N-l+1}\right), l = 2, \dots, N$, the {\it incomplete symmetrically dyadic block-sparsity pattern} of block size $k$, denoted by $\mathcal{IS}(N,l)$, is defined in accordance to the submatrix
$
    \mathbf \Sigma_{\mathcal I _{\mathcal{IS}}}$, { where }
$$
    \mathcal I _{\mathcal{IS}} = \mathcal{SD}(N)\setminus \bigcup_{l'=l}^{N} \bigcup_{r=1}^{2^{N-l'}}\mathbf J_{r,l'}.
$$

We observe that a matrix $\mathbf S\in \mathcal{IS}(N,l)$ is made of block‐diagonal symmetrically dyadic submatrices, i.e.,
$
    \mathbf S = \sum_{r=1}^{2^{N-l+1}}\mathbf S_{r}$,
where $\mathbf S_r$ is $\mathcal{SD}(l-1)$ in the $r$-th diagonal $k(2^{l-1}-1)\times k(2^{l-1}-1)$ block and zero otherwise.

Following the same logic, we define level-$l$ incomplete vertically dyadic block-sparsity pattern $\mathcal{IVD}(N, l)$ and its horizontal counterpart $\mathcal{IHD}(N, l)$ via
\begin{align}\label{eq:bspE}
    \mathcal I_{\mathcal{IVD}} =  \mathcal{VD}(N) \setminus \bigcup_{l^\prime = l} ^N \bigcup_{r=1}^{2^{N-l^\prime}} \mathbf I^\top_{r, l^\prime}, \,\,\,
    \mathcal I_{\mathcal{IHD}} = \mathcal{HD}(N) \setminus \bigcup_{l^\prime = l} ^N \bigcup_{r=1}^{2^{N-l^\prime}} \mathbf I_{r, l^\prime}.
\end{align}
As with $\mathcal VD(N)$ and $\mathcal HD(N)$, if a matrix $\mathbf V\in \mathcal{IVD}(N, l)$, then $\mathbf V^\top \in \mathcal{IHD}(N, l)$, and vice versa. Moreover, if matrices $\mathbf V\in \mathcal{IVD}(N, l)$ and $\mathbf H\in \mathcal{IHD}(N, l)$, then we can decompose them into block diagonal matrices
\begin{equation}
    \label{eq:decompVH}
    \mathbf V = \sum_{r=1}^{2^{N-l+1}}\mathbf V_r,\ \mathbf H = \sum_{r=1}^{2^{N-l+1}} \mathbf H_r,
\end{equation}
where the block-sparsity pattern of each $\mathbf{V}_r$ (respectively, $\mathbf{H}_r$) is given by $\mathcal{V D}(l-1)$ (resp. $\mathcal{H D}(l-1)$) in the $r$-th diagonal block of size $k\left(2^{l-1}-1\right)$ and zero elsewhere.

The following lemma summarizes the complexity of evaluating products of incomplete dyadic, elongated diagonal, and block diagonal matrices. They are crucial for establishing the main complexity result of this section (proof in the Appendix).

\begin{lemma}
    \label{lem:IHDED}
    Let $\mathbf H\in \mathcal{IHD}(N,l)$, $\mathbf V\in \mathcal{IVD}(N,l)$, $\mathbf D\in \mathcal{D}(N, l)$, and $\mathbf E\in \mathcal{ED}(N,l)$. Then $\mathbf H\mathbf E\in \mathcal{ED}(N,l)$ and $\mathbf V\mathbf E\in \mathcal{ED}(N,l)$, and the computation of either of these two products requires $O(ldk^2) $ flops, while $\mathbf E^\top \mathbf E\in \mathcal{D}(N,l)$ and $\mathbf E \mathbf D\in \mathcal{ED}(N,l)$, and the computation of either of these two products requires $O(dk^2)$ flops.
\end{lemma}

Next, we gather some basic facts about algebraic operations on matrices having dyadic structures and evaluate associated numerical complexities.

Suppose that $\mathbf H\in \mathcal{HD}(N)$, $\mathbf V\in \mathcal{VD}(N)$, and $\mathbf S\in \mathcal{SD}(N)$, then
\begin{equation*}
    \mathbf H = \sum_{l=1}^N \sum_{r=1}^{2^{N-l}} \mathbf H_{\mathbf J_{r,l}}, \,
    \mathbf V = \sum_{l=1}^N \sum_{r=1}^{2^{N-l}} \mathbf V_{\mathbf J_{r,l}}, \,
    \mathbf S = \sum_{l=1}^N \sum_{r=1}^{2^{N-l}} \mathbf S_{\mathbf J_{r,l}},
\end{equation*}
where the matrices in the sum have the same size as the left-hand matrices and have a block-sparsity pattern indicated by the index sets.

Let us recall that the computation of the product of two $k\times k$ matrices involves $O(k^3)$ flops, while adding them requires $O(k^2)$ flops.

\begin{proposition}
    \label{prop:SA}
    Let $\mathbf S\in \mathcal{SD}(N)$, and let $\mathbf A$ be a $d\times k$ matrix.
    Then
    \begin{equation*}
        \mathbf S \mathbf A=\sum_{l=1}^N \sum_{r=1}^{2^{N-l}} \mathbf S_{\mathbf J_{r,l}}\mathbf A_{I'_{r,l},\cdot}
    \end{equation*}
    and its computation involves $O\left(dk^2\log \left(d / {k}\right)\right)$ flops.
\end{proposition}
\begin{proof}
    The identity follows directly from the block-sparsity pattern. Thus, to obtain the computational cost, it is enough to find out the cost of multiplying symmetric cross-matrix $\mathbf S_{\mathbf J_{r, l}}$ by $\mathbf A_{I'_{r,l}}$. It involves the multiplication of two $k\times k$ matrices and then adding the results. In terms of complexity, we only need to consider the number of such multiplications, which is essentially the number of nonzero blocks in the cross-matrix $\mathbf S_{\mathbf J_{r,l}}$. There are $2^{l+1} - 3$ nonzero blocks in $\mathbf S_{\mathbf J_{r,l}}$. As a result, the total number of such matrix multiplications is
    \begin{equation}
        \sum_{l=1}^N\sum_{r=1}^{2^{N-l}} \left(2^{l+1} - 3\right)
        = (2N-3)2^N + 3
        = O\left(N2^N\right).
    \end{equation}
    Hence, the total number of flops required is
    \begin{equation}
        O\left(k^3N 2^N\right)
        = O\left(k^3\left(\frac{d}{k}+1\right)\log\left(\frac{d}{k}+1\right)\right) = O\left(dk^2\log \left(\frac{d}{k}\right)\right).
    \end{equation}
\end{proof}

\begin{proposition}
    \label{prop:VA}
    Let $\mathbf V\in \mathcal{VD}(N)$, and let $\mathbf A$ be a $d\times k$ matrix.
    Then
    \begin{equation*}
        \mathbf V \mathbf A=\sum_{l=1}^N \sum_{r=1}^{2^{N-l}} \mathbf V_{\mathbf I^\top_{r,l}}\mathbf A_{I'_{r,l},\cdot}
    \end{equation*}
    and its computation involves $O\left(dk^2\log \left({d}/{k}\right)\right)$ flops.
\end{proposition}

\begin{proof}
    The argument is almost identical to that for Proposition~\ref{prop:SA}, except that now we consider the column-matrix instead of cross-matrix. There are $2^l - 1$ nonzero blocks in $\mathbf V_{\mathbf I^\top_{r, l}}$, and the number of multiplications of two $k\times k$ matrices is $$
        \sum_{l=1}^N \sum_{r=1}^{2^{N-l}} (2^l-1)
        = N2^{N} - 2^N + 1  = O\left(N2^N\right).
    $$
    The rest is the same as the proof of Proposition \ref{prop:SA}.
\end{proof}

Analogously, we have
\begin{proposition}
    \label{prop:HA}
    Let $\mathbf H\in \mathcal{HD}(N)$, and let $\mathbf A$ be a $d\times k$ matrix. Then
    \begin{align*}
        \mathbf H \mathbf A & =\sum_{l=1}^N \sum_{r=1}^{2^{N-l}} \mathbf H_{\mathbf I_{r,l}}\mathbf A_{I'_{r,l},\cdot}
    \end{align*}
    and its computation involves $O\left(dk^2\log \left({d}/{k}\right)\right)$ flops.
\end{proposition}

The following result is tightly connected to the complexity of inverting dyadic matrices through Theorem~\ref{th:factor}.
\begin{theorem}
    \label{th:dya_mult}
    It takes $O(d^2k)$ flops to compute $\mathbf P \mathbf Q^\top$, where $\mathbf P,\mathbf Q\in \mathcal{VD}(N)$ and $d = k\left(2^{N}-1\right)$.
\end{theorem}

\begin{proof}
    We prove the result by mathematical induction, analyzing the recurrence for the number $a_N$ of flops needed for the multiplication.

    For $N=1$, the matrices $\mathbf P$ and $\mathbf Q$ are, in general, dense $k\times k$ matrices, so their product requires $a_1 = Ck^3$ flops for some constant $C>0$.

    For $N > 1$, the recursive relationship in iii) of Proposition~\ref{prop:propdyad} indicates the following sparsity patterns of products. Consequently, all blocks in $\mathbf P\mathbf Q^\top$ are potentially nonzero, and the corresponding products of blocks need to be evaluated
    \begin{equation*}
        \mathbf V(N+1)\mathbf V(N+1)^\top
        =  \begin{bmatrix}
            \mathbf V(N)\mathbf V(N)^\top +\mathbf 1_{d\times k} \mathbf 1_{k\times d} & \mathbf 1_{d\times k} \mathbf 1_{k\times k} & \mathbf 1_{d\times k} \mathbf 1_{k\times d}                               \\
            \mathbf 1_{k\times k}\mathbf 1_{k\times d}                                 & \mathbf 1_{k\times k}1_{k\times k}          & \mathbf 1_{k\times d}                                                     \\
            \mathbf 1_{d\times k} \mathbf 1_{k\times d}                                & \mathbf 1_{d\times k} \mathbf 1_{k\times k} & \mathbf V(N)\mathbf V(N)^\top+\mathbf 1_{d\times k} \mathbf 1_{k\times d}
        \end{bmatrix}.
    \end{equation*}
    It is thus easy to see the following recursive relationship
    \begin{equation*}
        a_{N+1} = 2a_N + C\left[4(2^N - 1)^2 + 4(2^N - 1) + 1\right]k^3= 2a_N + C\left(2^{N+1} - 1\right)^2k^3,
    \end{equation*}
    which leads to
    \begin{equation*}
        a_{N}
        =Ck^3 \sum_{K=1}^N 2^{N-K}\left(2^K - 1\right)^2  =Ck^3\left(2^{N}\left(2^{N+1}-2N-1\right)-1\right).
    \end{equation*}
    In particular, since $d = k(2^N-1)$, we have $a_N=O(d^2k)$.
\end{proof}


The following theorem presents the numerical complexity of our dyadic algorithm.

\begin{theorem}
    \label{th:maindyadic}
    Let $\mathbf \Sigma\in\mathcal{SD}(N)$. Algorithm~\ref{alg:dyadic} requires $O\left(dk^2\log^2 \left({d}/{k}\right)\right)$ flops to obtain the corresponding matrix $\mathbf P\in\mathcal{VD}(N)$.
\end{theorem}

\begin{proof}
    For $l = 2, \dots, N$, Algorithm~\ref{alg:dyadic} evaluates  matrices ${\mathbf P}^\top\check{\mathbf \Sigma}$, $\mathbf P \check{\boldsymbol \Sigma}^\prime$, $\boldsymbol \Sigma_{l,l} - {\check{\mathbf \Sigma}^\prime}{}^\top \check{\mathbf \Sigma}^\prime$, and $\mathbf A \mathbf G$.
    We analyze these computations one by one, repeatedly invoking Lemma~\ref{lem:IHDED}.

    First, for $\check{\mathbf \Sigma}^\prime = {\mathbf P}^\top\check{\mathbf \Sigma}$, we have $\mathbf P^\top \in \mathcal{IHD}(N, l)$ and $\check{\mathbf \Sigma}\in \mathcal{ED}(N, l)$. Hence, $\check{\mathbf \Sigma}^\prime\in \mathcal{ED}(N, l)$, and this operation requires $O(ldk^2)$ flops.
    Second, for
    $\mathbf A = \mathbf P \check{\mathbf \Sigma}^\prime$,
    we obtain $\mathbf A\in \mathcal{ED}(N, l)$ at a cost $O(ldk^2)$. Next, for
    $\tilde{\mathbf \Sigma} = \boldsymbol \Sigma_{l,l} - \check{\mathbf \Sigma}^\prime{}^\top \check{\mathbf \Sigma}^\prime $,
    as shown above
    this step incurs $O(dk^2)$ flops and
    $\tilde{\mathbf \Sigma} \in \mathcal{D}(N, l)$. Lastly, $\mathbf G$ is obtained by orthonormalizing $\tilde{\mathbf \Sigma}$ so that $\mathbf G\in \mathcal{D}(N, l)$, and it takes $O(dk^2)$ flops to calculate  $\mathbf A\mathbf G$. Therefore, for a given $l \ge 2$, the flop count is
    $   O(ldk^2) + O(ldk^2) + O(dk^2) + O(dk^2) = O(ldk^2)$.

    The Gram-Schmidt orthonormalization process costs $O(dk^2)$ flops, so the overall number of flops required by Algorithm~\ref{alg:dyadic} is
    \begin{equation*}
        \sum_{l=2}^N O(ldk^2)
        = O\left(d k^2 N^2\right)  = O\left(dk^2\log^2 \left({d}/{k}\right)\right).
    \end{equation*}

\end{proof}

Positive definite {\em block-tridiagonal matrices} with block size \(k\) form a notable subclass of \(\mathcal{SD}(N)\cap\mathcal{PD}(d)\). Such matrices have potentially nonzero \(k\times k\) blocks only on the main block diagonal and the first block super- and subdiagonals.
When $\mathbf \Sigma$ is block-tridiagonal, for each $l = 2,\dots,N$ the sparsity patterns of $\check{\mathbf \Sigma}$ and $\check{\mathbf \Sigma}^\prime$ can be refined. Define
$
    s(j,l) = (2^{l-1}-1)(2j-1),\quad 1\le j\le2^{N-l}.
$
Then, in iteration $l$,
$
    \widetilde{\mathcal{ED}}(N, l) = \{(i,j): 1\le j\le2^{N-l},\, s(j,l)\le i\le s(j,l)+1\}
$
describes the block-sparsity pattern of $\check{\mathbf \Sigma}$. Moreover, if we let
$$
    \mathcal T(l) = \{-2^{l-m}+1: m=2,\dots,l\} \cup \{2^{l-m}: m=2,\dots,l\},
$$
then it is shown in Lemma~\ref{lem:IHDED_band} that, in iteration $l$, the block-sparsity pattern of $\check{\mathbf \Sigma}^\prime$ is
$$
    \widetilde{\mathcal{ED}}^\prime(N, l) = \{(i,j): 1\le j\le2^{N-l},\, i=s(j,l)+t,\; t\in\mathcal T(l)\}.
$$
In this refined setting, dyadic algorithm complexity reduces by a factor of $\log(d/k)$.

\begin{lemma}
    \label{lem:IHDED_band}
    Let $\mathbf H\in \mathcal{IHD}(N,l)$, $\mathbf V\in \mathcal{IVD}(N,l)$, and $\mathbf E\in \widetilde {\mathcal{ED}}(N, l)$. Then $\mathbf H\mathbf E\in \widetilde {\mathcal{ED}}^\prime(N, l)$ and $\mathbf V\mathbf H\mathbf E\in \mathcal{ED}(N,l)$, with each product requiring $O(dk^2)$ flops.
\end{lemma}

\begin{theorem}
    \label{th:mainband}
    If $\mathbf \Sigma$ is block-tridiagonal with block size $k$, then Algorithm~\ref{alg:dyadic} requires $O\left(dk^2\log \left({d}/{k}\right)\right)$ flops to obtain the corresponding matrix $\mathbf P\in\mathcal{VD}(N)$.
\end{theorem}
\begin{proof}
    By Lemma~\ref{lem:IHDED_band}, for $l = 2, \dots N$, $O(dk^2)$ flops are required to calculate $\mathbf P^\top \check{\mathbf \Sigma}$ and $\mathbf P\check{\mathbf \Sigma}^\prime$. As a result, the total number of flops required by Algorithm~\ref{alg:dyadic} is
    \begin{equation*}
        \begin{split}
            \sum_{l=2}^N O(dk^2)
             & = O(dNk^2) =O \left(dk^2\log\left({d}/{k}\right)\right)
        \end{split}
    \end{equation*}
\end{proof}

\section{Packing sparse matrices}
\label{sec:packing}
We present a method for packing sparse matrices into block-tridiagonal matrices and investigate its use in recursively recovering the dyadic structure.
The R implementation of this algorithm can be found in the GitHub repository \cite{packing_algorithm_github}.
The main focus is thus on finding, for a given PD $d\times d$ matrix $\boldsymbol \Sigma$, a permutation matrix $\boldsymbol\pi$ such that $\boldsymbol\pi^\top \boldsymbol \Sigma \boldsymbol\pi$ is block-tridiagonal, allowing the application of the dyadic factorization algorithm.
However, since the width of the sought dyadic (in this case block-tridiagonal) matrix is not known \textit{a priori} and the developed method uses a criterion that promotes tight packing along the diagonal, the algorithm applies to more general incomplete dyadically structured matrices.

From now on, only the zero–nonzero pattern of $\boldsymbol \Sigma$ matters; accordingly, we treat $\boldsymbol \Sigma$ as a symmetric $0$-$1$ matrix.
Moreover, for a permutation $\pi$, the matrix representing it is denoted by $\boldsymbol \pi$, i.e.,
$
    \boldsymbol \pi =[\delta_{j,\pi(i)}]_{i,j=1}^d,
$
where $\delta_{kj} = 1$ if $k = j$ and zero otherwise.
Trivially, any $d \times d$ matrix is a block-tridiagonal matrix with one $d \times d$ block. The objective is to find a $\pi$ that packs the nonzero entries as close to the diagonal as possible using a natural criterion.

\subsection{Optimal permutation}
\label{subsec:optper}
The essence of the dyadic algorithm is to optimally exploit the knowledge of zero blocks located at specific distances from the diagonal of the matrix $\boldsymbol\Sigma$. In the following definition, we introduce a criterion for permutation optimality that favors permutations placing nonzero elements near the diagonal.

\begin{definition}
    \label{def:band}
    For a symmetric $0$-$1$ matrix $\boldsymbol\Sigma$, a permutation $\pi$, and row indices $i$, $j$, we call
    \begin{itemize}
        \item $D_i = \{j: \boldsymbol{\Sigma}_{i,j} \neq 0 \}$ --- the $i$-th neighborhood, i.e. the set of neighbors of row $i$,
        \item $l^{(\pi)}_i = \max\{|\pi(i)-\pi(j)|:j \in D_i\}$ ---
              the half-width of $D_i$ under $\pi$,
        \item
              $\Vert \boldsymbol l^{(\pi)}\Vert_1 = \sum_{i=1}^d l^{(\pi)}_i$ --- the half-width of $\boldsymbol\Sigma$ under $\pi$,
        \item
              a permutation $\pi$ is called \textbf{optimal} if it belongs to the following set
              \begin{equation}\label{crit 1}
                  \Pi(\boldsymbol{\Sigma}) = \underset{\pi }{\operatorname{argmin}} \Vert \boldsymbol l^{(\pi)}\Vert_1.
              \end{equation}
    \end{itemize}
\end{definition}

Since exactly solving \eqref{crit 1} is computationally infeasible for large $d$, our goal is to approximate the optimal permutation. Our proposed approach leverages the relationship between the optimal permutation and its associated distance matrix.

\subsection{Distance matrices}
\label{sec: dist_matrix}
Distance matrices are widely used in modern multivariate data analysis and serve as the foundation of our approach. A one-dimensional distance matrix determines its configuration up to translation and reflection. In particular, for a permutation $\pi$, the matrix
\[
    \mathbf G_\pi
    =
    [|\pi(i)-\pi(j)|]_{i,j=1}^d
\]
determines the corresponding ordering up to reversal. Algorithm~\ref{alg:config} reconstructs this ordering, while the general reconstruction result is provided in Section~\ref{DistMat_supp} of the Appendix.


\begin{algorithm}
    \caption{Extracting configuration from the permutation distance matrix}
    \label{alg:config}
    \begin{algorithmic}\vspace{2mm}

        \STATE{{\sc Input:} $\mathbf G_{\pi} = [G_{ij}]$ - $d\times d$ distance matrix associated with a permutation $\pi$}\vspace{1mm}
        \STATE{$\alpha(1) =\min \{i:\, \#\{j: G_{ij}=1\}=1\}$};  \COMMENT{{\tt \small Two i's satisfy the condition}}\vspace{1mm}
        \STATE{$\alpha(2) = find(i: G_{i, \alpha(1)}=1)$}; \COMMENT{{\tt \small One i satisfies the condition}}\vspace{1mm}
        \STATE{ $m=2$};
        \WHILE{$ m < d$}\vspace{2mm}
        \STATE{$\alpha(m+1) = find(i\neq \alpha(m-1): G_{i, \alpha(m)}=1)$}; \COMMENT{{\tt \small One i satisfies the condition}}\vspace{1mm}
        \STATE{$m=m+1$}\vspace{1mm}
        \ENDWHILE \vspace{1mm}
        \STATE{{\sc Return:} $\alpha$ - a permutation whose distance matrix is $\mathbf G_\pi$ }
    \end{algorithmic}
\end{algorithm}

\begin{remark}
    \label{rem:1}
    For the permutation distance matrix $\mathbf G_\pi$, all nearest-neighbor distances equal one. Hence, Algorithm~\ref{alg:config} recovers the ordering from $2d-2$ nearest-neighbor relations, up to reversal. If $\rho(i)=d+1-i$, then the recovered ordering is either $\pi$ or $\rho\circ\pi$. Since reversal preserves every half-width, $\pi$ is optimal if and only if $\rho\circ\pi$ is optimal.
\end{remark}

\subsection{The approximate distance matrix optimization method}
We outline the core ideas of our method, which consist of two steps: first, approximating the optimal permutation distance matrix with the ``distance matrix'' \(\mathbf A_{\boldsymbol \Sigma}\) and second, approximating the optimal permutation using \(\mathbf A_{\boldsymbol \Sigma}\). Since the optimal permutation distance matrix $\mathbf G_\pi$ is unknown, we approximate it using only the row neighborhoods. Because the resulting distances need not be Euclidean, we apply one-dimensional classical multidimensional scaling (CMDS) \cite{Torgerson52,GOWER66} and rank the resulting coordinates to estimate the permutation. A general review of CMDS is provided in the Appendix. For a symmetric $0$-$1$ matrix $\boldsymbol\Sigma$, the ``distance matrix'' is given by

\begin{equation}
    \label{SymmDiffDist}
    \mathbf A_{\boldsymbol \Sigma} = \left[ \frac{\left |D_i \triangle D_j\right |}{2}\right]_{i,j=1}^d,
\end{equation}
where $\left |D_i \triangle D_j\right |$ is the number of elements in the symmetric difference of sets, which is a metric on the space of all subsets of $\{1,\dots, d\}$. We show that, under certain conditions, $\mathbf A_{\boldsymbol \Sigma}$ serves as a good approximation of $\mathbf G_\pi$ for the optimal $\pi$.
The approximation error is bounded in Theorem~\ref{T:1.3} and the subsequent conclusions.

A simple measure of the proximity of two matrices is the average absolute difference between the corresponding entries. However, as seen in the previous section, determining the permutation $\pi$ does not require knowledge of all $\mathbf G_\pi$ elements. It is sufficient to know the nearest neighbor on each side of the point corresponding to a given row (see Remark~\ref{rem:1}).
Therefore, we consider the following criterion
\begin{equation}
    \label{matrix_dist}\delta_m(\mathbf {G}_\pi,\mathbf A_{\boldsymbol \Sigma})
    = \frac{1}{dm}
    \sum_{\substack{i,j=1 \vspace{0.5mm}\\
            \left[\mathbf G_\pi\right]_{ij} \le m}}^d
    \left |\left[\mathbf G_\pi\right]_{ij} - \left[\mathbf A_{\boldsymbol \Sigma}\right]_{ij}\right |.
\end{equation}

The constraint $[\mathbf G_\pi]_{ij} \le m$ means that we consider only the nearest neighbors at a distance no greater than $m$.
We are interested in determining the set of permutations
\begin{equation}\label{crit 2}
    B_{\boldsymbol \Sigma, m} = \underset{\pi }{\operatorname{argmin}} \;\;\delta_m(\mathbf{G}_\pi,\mathbf A_{\boldsymbol \Sigma}).
\end{equation}
Ideally, we would like to know whether the optimal permutation is in $B_{\boldsymbol \Sigma, m}$, which will make $\mathbf A_{\boldsymbol \Sigma}$ a good candidate for approximation.
Such a strong result is unknown, but we instead derive an upper bound for $\delta_m(\mathbf{G}_\pi,\mathbf A_{\boldsymbol \Sigma})$ which attains its minimum for the optimal $\pi$ (see Remarks~\ref{rem: 2b}  and~\ref{rem:2}).
In the following theorem and proposition, we use the notation of Subsection~\ref{subsec:optper}. The proof can be found in the Appendix.
\begin{theorem} \label{T:1.3}
    Let $\boldsymbol{\Sigma}$ be a symmetric $0$-$1$ matrix and let $\pi$ be an arbitrary permutation. Suppose $i$ and $j$ are two distinct row indices such that
    \begin{equation}
        \label{eq:ass}
        |l^{(\pi)}_i - l^{(\pi)}_j|\le \left[\mathbf  G_\pi\right]_{ij}\le l^{(\pi)}_i + l^{(\pi)}_j.
    \end{equation}
    Then
    \begin{equation}
        \left |\left[\mathbf G_\pi\right]_{ij} - \left[\mathbf A_{\boldsymbol \Sigma}\right]_{ij}\right | \le \frac{(2l^{(\pi)}_i+1-|D_i|) +(2l^{(\pi)}_j+1- |D_j|)}2 .
        \label{eq: approx2}
    \end{equation}
\end{theorem}

It is easy to see that Theorem \ref{T:1.3} allows us to describe how close the matrices $\mathbf A_{\boldsymbol \Sigma}$ and $\mathbf G_\pi$ are for an arbitrary permutation $\pi$ (compare \ref{matrix_dist} with \ref{eq: approx2}). However, before discussing its consequences in this context, we focus on the assumptions and conclusions of the theorem.
Since assumptions \eqref{eq:ass} involve the unknown  $l^{(\pi)}_i$ and $l^{(\pi)}_j$, it is practical to establish an alternative condition that directly implies these relations. The following proposition provides a sufficient condition for the assumption on the right-hand side.

\begin{proposition} \label{prop: assumpt.1}
    For a permutation $\pi$ and indices $i$ and $j$:
    $$
        D_i\cap D_j \neq \emptyset
        \Rightarrow \left[\mathbf G_\pi\right]_{ij}\le l^{(\pi)}_i + l^{(\pi)}_j.
    $$
\end{proposition}
\begin{proof}
    Let $t\in D_i\cap D_j \neq \emptyset$.
    By the definition of half-width
    $$
        \left[\mathbf G_\pi\right]_{ij} = |\pi(i) - \pi(j)| \le |\pi(i) - \pi(t)| + |\pi(t) - \pi(j)| \le l^{(\pi)}_i + l^{(\pi)}_j.
    $$
\end{proof}

The advantage of condition $D_i\cap D_j \neq \emptyset$ is that it is easy to check and does not refer to the optimal permutation. Furthermore, the condition has a clear interpretation (rows $i$ and $j$ must have a common neighbor). However, this condition is rather restrictive, as only relatively close rows can satisfy it (in the following section, we discuss how this limitation can be circumvented).
When the condition is not satisfied, $ \left[\mathbf A_{\boldsymbol \Sigma}\right]_{ij} $ can be significantly smaller than $ [\mathbf G_\pi]_{ij} $, as illustrated in the following example.

\begin{example}
    \label{ex:tridiagonal}
    For full symmetric tridiagonal matrices, the set of optimal permutations consists of two elements: the identity permutation and the permutation $\rho$, defined as  $\rho(i)=d+1-i$ for $i=1,2,\dots,d$.
    Moreover, the following relations hold
    \begin{itemize} \label{ex: tridiag}
        \item[\it i)] For any pair $i,j$, $ \left[\mathbf A_{\boldsymbol \Sigma}\right]_{ij} \le 3$,   (proof: $|D_i \triangle D_j|/2  \le (|D_i| + |D_j|)/2 \le 3$)
              .
        \item[\it ii)]  For an optimal  $\pi$, if $\pi(i),\pi(j) \notin \{1,d\}$ and $D_i\cap D_j \neq \emptyset$, then $[\mathbf G_\pi]_{ij} = \left[\mathbf A_{\boldsymbol \Sigma}\right]_{ij}$.
    \end{itemize}

    From {\it i)}, even for moderate $d$, many index pairs $i, j$ satisfy $[\mathbf G_\pi]_{ij} >  \left[\mathbf A_{\boldsymbol \Sigma}\right]_{ij}$.
    On the other hand, {\it ii)} shows that the approximation is exact for the closest neighbors.
\end{example}

\begin{remark}
    \label{rem:neigh}
    Mathematically, condition $D_i\cap D_j \neq \emptyset$ is sufficient for the assumption on the right side to hold. However, in practice, relying on the approximation when the condition is not fulfilled may be unreliable, as it may result in a significant underestimation of $[\boldsymbol G_\pi]_{ij}$.
    Notably, the condition is guaranteed for elements $i$ and $j$ of any nearest neighbors set $ D_r$, since $ i,j \in D_r $ implies $ r \in D_i \cap D_j $. Furthermore, information about the nearest neighbors is sufficient (see Remark~\ref{rem:1}).
\end{remark}

There is no simple alternative to the assumed left-hand side bound in \eqref{eq:ass}. Nevertheless, if the optimal permutation leads to a dense {block-tridiagonal} matrix, it is reasonable to expect that the parameters $ l^{(\pi)}_i $ and $ l^{(\pi)}_j $ are similar or relatively small. Thus, the left-hand side is expected to be small.

The upper bound \eqref{eq: approx2} averages two similar terms. To better understand the factors influencing this bound, we introduce the concept of density for the sets involved.
From the definition of $l_i^{(\pi)}$, it follows that
$
    \pi(D_i) \subseteq \pi(i) + C(l^{(\pi)}_i)$,
where $C(t) = \{-t, \dots, -1, 0, 1, \dots, t\}$. The cardinalities of these sets are $|D_i|$ and $2l^{(\pi)}_i +1$, respectively, with the obvious inequality $|D_i| \le 2l^{(\pi)}_i +1$.
To quantify how densely the set $\pi(D_i)$ is distributed within $\pi(i) + C(l^{(\pi)}_i)$, we define the density measure
$\eta_i(\pi) = \frac{|D_i|}{2l^{(\pi)}_i +1}$.
This quantity satisfies $\frac{1}{2l^{(\pi)}_i +1} \le \eta_i(\pi) \le 1$, where values closer to 1 indicate a higher density of elements from $\pi(D_i)$ in the given range. This allows density comparisons across different indices on the same scale.
To assess the overall density across all rows, we define
$\bar{\eta}(\pi) = \frac{1}{d} \sum_{i=1}^{d} \eta_i(\pi)$.

\begin{remark}
    \label{rem:imp}
    The upper bound \eqref{eq: approx2} decreases as the densities $\eta_i(\pi)$ and $\eta_j(\pi)$ increase. In other words, the denser $\pi(D_i)$ and $\pi(D_j)$ are within their respective supersets $\pi(i) + C(l^{(\pi)}_i)$ and $\pi(j) + C(l^{(\pi)}_j)$, the closer $\left[\mathbf A_{\boldsymbol \Sigma}\right]_{ij}$ and $\left[ \mathbf G_\pi\right]_{ij}$ become. In the extreme case where $\eta_i(\pi) = 1$ and $\eta_j(\pi) = 1$, we obtain $
        \left[\mathbf A_{\boldsymbol \Sigma}\right]_{ij} = \left[\mathbf G_\pi\right]_{ij}$.
    Thus, the density measure $\eta_i(\pi)$ enables consistent comparisons across different indices and is directly related to the accuracy of the approximation \eqref{eq: approx2}.
\end{remark}

The repeated use of inequalities \eqref{eq: approx2} applied to the formula \eqref{matrix_dist} for
$\delta_m(\mathbf {G}_\pi,\mathbf A_{\boldsymbol{\Sigma}})$
leads to the following important proposition:
\begin{proposition} \label{Prop:2}
    Let  $\pi$ be an arbitrary permutation.
    Suppose that \eqref{eq:ass} is satisfied for any pair of row indices $i,j \in \{1,...,d\} $ such that $[\mathbf {G}_\pi]_{ij}\le m$ for some natural number $m$.
    Then the following inequality holds
    \begin{equation}
        \begin{aligned}
            \delta_m(\mathbf {G}_\pi,\mathbf A_{\boldsymbol{\Sigma}})
            \le
            \frac{1}{d}\left\{2Q(\pi,\boldsymbol \Sigma)
            - \left[
            \sum_{\pi(i)=1}^{m}\frac{m+1-\pi(i)}{m} \left(2l^{(\pi)}_{i}+1-|D_{i}|\right) \right.\right. \\
                \left.\left.
                +
                \sum_{\pi(i)=d-m+1}^d  \frac{m-d+\pi(i)}{m}\left(2l^{(\pi)}_{i}+1-|D_{i}|\right)\right]\right\},
            \label{eq:uperr bound 2}
        \end{aligned}
    \end{equation}
    where $Q(\pi,\boldsymbol \Sigma) =2\Vert \boldsymbol l^{(\pi)}\Vert_1 +  \sum_{i=1}^d(1-|D_i|)$. Furthermore, when $m=d$, then
    \begin{equation}
        \delta_d(\mathbf {G}_\pi,\mathbf A_{\boldsymbol{\Sigma}}) \le \frac{(d-1)}{d^2}Q(\pi,\boldsymbol \Sigma)
        \label{eq:uperr bound 1}
    \end{equation}
\end{proposition}
The proof of this proposition is in the Appendix.

\begin{remark} \label{rem: 2b}
    Let us first focus on the case where $m = d$, since in this setting the upper bound formula (\ref{eq:uperr bound 1}) is simpler and easier to interpret.   We highlight two key observations.
    Firstly, $Q(\pi,\boldsymbol \Sigma)$ consists of two components. The first is the scaled half-width of the $\boldsymbol{\Sigma}$ matrix under the permutation $\pi$ ($2\Vert \boldsymbol l^{(\pi)}\Vert_1$), while the second component ($\sum_{i=1}^d(1-|D_i|)$) is independent of the choice of permutation. Consequently,
    $\underset{\pi }{\operatorname{argmin}} \;\;Q(\pi,\boldsymbol \Sigma) = \underset{\pi }{\operatorname{argmin}} \;\;\Vert \boldsymbol l^{(\pi)}\Vert_1.
    $
    This equation shows that the distance $\delta_d(\mathbf {G}_\pi,\mathbf  A_{\mathbf \Sigma})$ has the tightest upper bound for the optimal permutation.
    Secondly, the upper bound by the definition is the sum of bounds \eqref{eq: approx2} over all pairs of indices $i\ne j$. In light of Remark \ref{rem:imp}, if for the permutation $\pi$, the mean $\bar \eta(\pi)$ is close to 1, then the upper bound should be relatively close to 0.

    Given the above observations,  $\mathbf {A_{\boldsymbol \Sigma}}$ is a good candidate for the approximation of $\mathbf {G}_\pi$ associated with the optimal $\pi$, at least when certain conditions are met.
\end{remark}

When $m=d$, we assume that \eqref{eq:ass} holds for all index pairs. Example \ref{ex:tridiagonal} shows that this assumption is quite restrictive.
The next remark concerns the properties of the upper bound formula \eqref{eq:uperr bound 2} for small $m$, meaning that condition \eqref{eq:ass} is satisfied only for the nearest neighbors — a setting that is of particular interest from our perspective.

\begin{remark}\label{rem:2}
    Even though the formula for the upper bound given in \eqref{eq:uperr bound 2} is more complex than \eqref{eq:uperr bound 1}, when $ m/d $ is small, the negative term is negligible. Thus, it is reasonable to expect that the upper bound is the smallest for a permutation close to the optimal permutation, as its properties mainly depend on $Q(\pi,\boldsymbol \Sigma)$  (see observation in Remark \ref{rem: 2b}). This suggests that the matrices $\mathbf A_{\boldsymbol \Sigma}$ and $\mathbf G_\pi$ associated with the optimal permutation may be close in the metric $\delta_m$ when $m$ is relatively small compared to $d$.
\end{remark}
\begin{remark}
    The matrix proximity condition is primarily theoretical. Simulations show that even when $\bar\eta(\pi)$ deviates significantly from 1, our method, which relies on the approximation matrix $\mathbf{A}_{\boldsymbol{\Sigma}}$, still yields near-optimal permutations for a broad class of matrices, suggesting robustness to violations of this condition.
\end{remark}

\subsection{Algorithm}
Remark \ref{rem:1} suggests that knowing the distances between the nearest neighbors is sufficient to determine the configuration of points. Motivated by this, we present an algorithm that computes the configuration of points based solely on the distances between their nearest neighbors.

The algorithm consists of three steps. We extract indices on which the permutation will be built ({\sc Skeleton}), attach to them the local configurations ({\sc Flesh}), and use these to build the global configuration ({\sc Body}).
A more detailed breakdown of these components is provided below. \vspace{1mm}

\begin{description}
    \item[\sc Skeleton:]  A subset of strategically chosen row indices $J = \{ i_1, i_2, \dots, i_{|J|} \}$, selected based on a criterion discussed later.
    \item[\sc Flesh:]  A collection of local configurations, where for each $ i \in J $, CMDS determines the arrangement of points in $D_i$ using the corresponding submatrix of $\mathbf{A} = \mathbf{A}_{\boldsymbol{\Sigma}} $. These local configurations are then used to construct the full permutation.
    \item[\sc Body:]
          A global configuration and the associated permutation, assembled from local configurations. Since the local configurations are determined only up to translation and reflection, alignment along the skeleton $J$ is required.\vspace{1mm}
\end{description}

\subsubsection{Local configurations} \label{sec: loc_config}
The selection of elements in $J$ is closely tied to the alignment process used for the final configuration. Therefore, we first focus on local configurations ({\sc Flesh}).
Here, we apply CMDS on sets $D_i$, for all $i \in J$. More precisely, let $D_i = \{m_1,...,m_{s_i}\}$, where $m_k<m_{k+1}$.
Based on the neighborhoods of the $D_i$ elements, we calculate the distance matrices
$\left[\mathbf A^{(i)}\right ]_{kj} = \frac{|D_{m_k} \triangle D_{m_j}|}2,\,\,\,  k,j \in \{1,...,s_i\}
$.
Next, we apply CMDS to obtain a one-dimensional configuration vector $\mathbf x^{(i)} \in \mathbb{R}^{s_i}$, where $x^{(i)}_k$ represents the position of row index $m_k \in D_i$ on a line, as determined by CMDS.
The reason for using $\mathbf A^{(i)}$ is to effectively use the approximation \eqref{eq: approx2}. As discussed, it is beneficial that $ D_{m_k} \cap D_{m_j} \neq \emptyset $ and, by Remark~\ref{rem:neigh}, this condition is valid for any pair of $ D_i $ elements.
We thus treat the sets $ D_i $, $i\in J$ as fundamental building blocks where CMDS provides a reasonable configuration of points.

\subsubsection{Global configuration}
Each vector $\mathbf x^{(i)}$ contains only partial information about the positions of the indices, since it captures only relative relationships within $D_i$. To obtain a global configuration from this ``flesh", we unify information in all $\mathbf x^{(i)}$ utilizing the ``skeleton" $J$.

We define the $d\times |J|$ matrix $ \mathbf Y=\left[y_{ij}\right] =\left[\mathbf y_{\cdot j}\right] $, which encodes local configurations for indices $ i_1, i_2, \dots, i_{|J|} \in J $. Each entry is given by
$$
    y_{m_kj} = \begin{cases}
        x^{(i_j)}_k, & \text{ if }  m_k \in D_{i_j} \\
        0,           & \text { otherwise},
    \end{cases}
$$ where $x^{(i_j)}_k$ is the position on a line of a row index $ m_k \in D_{i_j} $.
The alignment process iteratively transforms subsequent columns $\left[\mathbf y_{\cdot  j} \right]$ of the matrix $\mathbf Y$ and stores the results of iterations in a matrix $\widetilde{\mathbf Y}$.
\begin{algorithm}
    \caption{Flesh-to-body}
    \label{alg:flesh2body}
    \begin{algorithmic}\vspace{2mm}
        \STATE{{\sc Input:} $\mathbf Y =\left [y_{ij}\right]=\left[\mathbf y_{\cdot j}\right]$ - $d\times |J|$  local configurations (columnwise) along the skeleton,
        $D_{i_1}, D_{i_2},...,D_{i_{|J|}}$ nearest neighbors sets of indices in the skeleton} \vspace{1mm}
        \STATE{$ \widetilde{\mathbf{Y}}_{1} = \left[ \mathbf y_{\cdot 1} \right] $};  \COMMENT{{\tt \small Initialization of a new configuration matrix}}\vspace{2mm}
        \FOR{$k=2$ to $|J|$}\vspace{2mm}
        \STATE{$(\hat a, \hat b) = \underset{(a,b)\in \{-1,1\}\times \mathbb R}{\operatorname{argmin}} f\left(a,b; \mathbf y_{\cdot k}, \widetilde{\mathbf Y}_{k-1}\right)$};
        \COMMENT{{\tt \small The objective function in \eqref{eq:rec_align}}}
        \vspace{1mm}
        \FOR{$j=1$ to $d$}\vspace{2mm}
        \STATE{$ \widetilde{y}_{jk} = (\hat a  y_{j k} + \hat b) \mathds{1}(j \in D_{i_k})$};
        \COMMENT{{\tt \small  Adding the k-th column}} \vspace{1mm}
        \ENDFOR \vspace{1mm}
        \ENDFOR \vspace{1mm}
        \STATE{${\mathbf x}= {\rm repeat}({\tt NA},d)$};  \COMMENT{{\tt \small Initialization of the final configuration vector}}\vspace{2mm}
            \STATE{${\mathbf 1}= {\rm repeat}(1,|J|)$};  \COMMENT{{\tt \small Column vector of ones}}\vspace{2mm}
        \FOR{$j=1$ to $d$}  \vspace{1mm}
        \STATE{$n_j =  \left|\{k;j\in D_{i_k}\}\right|$};  \vspace{1mm}
        \IF {$n_j>0$}\vspace{1mm}
        \STATE{$ x_{j}= \widetilde {\mathbf y}_{j \cdot} \mathbf 1/n_j$}; \COMMENT{{\tt \small Averaging available local coordinates}}\vspace{1mm}
        \ENDIF\vspace{1mm}
        \ENDFOR
        \vspace{3mm}
        \STATE{{\sc Return:} $ {\mathbf x}=(x_j)_{j=1}^d $ - The global configuration along the skeleton based on $\mathbf Y$.  }
    \end{algorithmic}
\end{algorithm}
In the algorithm, the first configuration vector $\widetilde{\mathbf y}_{\cdot 1}$ equals $\mathbf y_{\cdot 1}$.
Then the subsequent configuration $\widetilde{\mathbf y}_{\cdot k}$ given the previous $\widetilde{\mathbf Y}_{k-1}=\left[\widetilde{\mathbf y}_{\cdot i}\right]_{i<k}$ is obtained by minimizing in $(a,b)\in \{-1,1\}\times \mathbb R$ the following objective function
\begin{equation}
    \label{eq:rec_align}
    f(a,b; \mathbf y_{\cdot k}, \widetilde{\mathbf Y}_{k-1})
    = \sum_{m = 1}^{k-1} \sum_{z = 1}^{d} (\tilde y_{z m} - a y_{z k} - b)^2 \mathds{1}(z \in D_{i_m} \cap D_{i_k})
\end{equation}

As seen in Algorithm~\ref{alg:flesh2body}, the unified configuration vector $\mathbf x$ is obtained by row-wise averaging local configurations given by $\widetilde{\mathbf Y}$, producing a sequentially aligned global configuration.
For $ k \in \bigcup_{j = 1}^{|J|} D_{i_j}$, the entries $x_k$ of the configuration vector $\mathbf x$ can be written as
$
    x_k = \frac 1{m_k}
    \sum_{j = 1}^{|J|}
    \tilde y_{kj} \mathds{1}(k \in D_{i_j})
$, where $m_k$ is the number of $D_{i_j}$'s, $j=1,\dots, |J|$, containing $k$.

If the index $ k $ does not belong to any of the sets $  D_{i_j} $, the corresponding element in the configuration $ \mathbf x $ is not available, which should be prevented by a proper choice of the skeleton $J$. Below, we describe the conditions that the set $ J $ must satisfy to avoid this problem.
If $\mathbf{x}$ is fully determined, we approximate the optimal permutation as $\hat \pi=(r_1(\mathbf{x}), \dots , r_d(\mathbf{x}))$, where $r_i(\mathbf{x})$ denotes the position of $x_i$ in the ordered sequence.

\begin{remark}
    Let us examine the objective function given in \eqref{eq:rec_align}.
    This function measures the alignment of a new vector with previously modified vectors after reflection $a$ and translation $b$, but only over the indices common to $D_{i_k}$ and $D_{i_m}$. The optimization procedure seeks the transformation that best fits prior data in the least-squares sense.
    It is closely related to the Procrustes transformation, which, in general, aligns two point sets in Euclidean space using translation, scaling, and rotation to minimize the sum of squared distances between corresponding points. See \cite{BorgGroenen05} for details.
    For $k=2$, our procedure is a special case, involving only translation and reflection in one dimension. For $k>2$, it extends this idea to multiple vectors.

    By solving the least-squares optimization, $\hat a$ is the sign of
    \begin{equation}
        \label{eq:optimal_ab}
        \sum_{j=1}^{k-1} \sum_{s \in D_{i_j} \cap D_{i_k}}
        \left( n_k\tilde y_{s j} - \sum_{m=1}^{k-1} \sum_{z \in D_{i_m} \cap D_{i_k}} \tilde  y_{z m}\right)
        \left( n_ky_{s k} - \sum_{m=1}^{k-1} \sum_{z \in D_{i_m} \cap D_{i_k}} y_{zk}\right)
    \end{equation}
    while
    $
        \hat{b} = \sum_{m=1}^{k-1} \sum_{z \in D_{i_m} \cap D_{i_k}} (\tilde y_{zm} - \hat a y_{zk})/n_k
    $,
    where $n_k=\sum_{m=1}^{k-1} \left| D_{i_m} \cap D_{i_k}\right|$.
\end{remark}
\subsubsection{The choice of skeleton}
The choice of skeleton $J$ is crucial and serves as the first step of the packing algorithm.
The main challenge is to find a skeleton that has sufficient overlaps between the neighborhoods of its elements when the matrix is sparse.
Our method does not assume any specific matrix structure, except that the graph induced by the sparsity pattern is connected. This entails no loss of generality, since disconnected components can be treated separately.
We select the skeleton $J$ that satisfies the following key criteria. \vspace{1 mm}
\begin{description}
    \item[\it Complete skeleton neighborhoods:] The skeleton neighborhoods cover all indices, i.e.,
          $
              \bigcup_{j\in J} D_j = \{1, 2, \ldots, d\}
          $.
    \item[\it Robust estimators:] The estimator $\hat a$ calculated at each step of alignment is robust, i.e., there is sufficient overlap between the skeleton neighborhoods.
    \item[\it Computational complexity and accuracy:]
          The calculation of the estimators should be efficient, balancing complexity and accuracy. The set $J$ should be as small as possible to reduce computational cost.\vspace{1 mm}
\end{description}

The first criterion must be satisfied as failing to do so would leave some points without assigned positions in the global configuration vector, making the final permutation incomplete.
Increasing the number of elements in the skeleton always ensures that this criterion is met, even for very sparse matrices. However, the challenge is not just to satisfy this criterion, but also to ensure sufficient overlap between the neighborhoods of skeleton elements, which is associated with the second criterion.

The second criterion is critical for achieving a high-quality approximation to the optimal permutation. Simulations indicate that an incorrect estimation of the parameter $ a $ can result in a permutation $ \pi $ where the half-width of $ \boldsymbol{\Sigma} $ under $ \pi $ is significantly larger than in the optimal case. See Example~\ref{ex:cont} for further discussion.
In this context, an analysis of the objective function of the alignment process \eqref{eq:rec_align} shows that for every $r \geq 2$, a necessary condition to uniquely determine the estimators is
\begin{equation} \label{est_cond}
    \sum_{m = 1}^{r-1} |D_{i_m} \cap D_{i_r}| \geq 2.
\end{equation}
The necessity of \eqref{est_cond} arises from the fact that two parameters need to be estimated. If \eqref{est_cond} is not satisfied, only one observation is available, making a unique estimate impossible. Moreover, a stronger condition enhances robustness: for each $r$, having at least two different indices in $D_{i_r} \cap \bigcup_{m = 1}^{r-1} D_{i_m}$ leads to more stable estimators in the alignment step.  This stronger condition holds when $i_r$ is a direct neighbor of a previous skeleton element $i_m$ ($m<r$), since ${i_r,i_m} \in D_{i_m} \cap D_{i_r}$.
However, if the elements of the skeleton are not direct neighbors, they may only be connected through a common neighbor outside the skeleton, in which case the sum in \eqref{est_cond} may be equal to one. This occurs in Example~\ref{ex: tridiag} with a full tridiagonal matrix. In this particular case, if $i_1=1$ and $i_2=3$, then $|D_{i_1}\cap D_{i_2}|=1$.

Based on the above discussion, it is desirable that successive elements of the sequence $ i_r $ are elements of $\bigcup_{m = 1}^{r-1} D_{i_m}$.
The above condition has another positive aspect. For a given $ m $, if the indices $ i_r $ and $ i_m $ are nearest neighbors, it can be expected that $ D_{i_r} \cap D_{i_m} $ contains more elements than the pair $\{i_r,i_m\}$. Evidently, this improves the precision of parameter estimation in the alignment process.
\begin{example}[Cont.]
    \label{ex:cont}
    Consider two scenarios based on Example~\ref{ex:tridiagonal} with different choices of the skeleton sequence $J = {i_1, i_2, \dots, i_{|J|}}$.

    In the first scenario, $ i_j = j+1 $ so $ D_{i_j} = \{i_j-1, i_j, i_j+1\} $, for $ j = 1, \dots, d-2 $. With this choice of $ J $, for $ r > 1 $, $ \sum_{m=1}^{r-1} |D_{i_m} \cap D_{i_r}| = 2 $. Since the neighborhoods $ D_{i_1}, \dots, D_{i_{r-2}} $ are disjoint from $ D_{i_r} $, they do not affect parameter estimation. In addition, $ D_{i_r} $ shares exactly two neighbors (indices $ i_r-1 $ and $ i_r $) with neighborhood $ D_{i_{r-1}} $.
    Consequently, in the alignment process, the parameter $ a $ is selected at each step in such a way that the positions of the indices $ i_r-1 $ and $ i_r $ in local configurations $ x^{(r-1)} $ and $ x^{(r)} $ maintain the same relative ordering. Furthermore, it can be shown that, in practical terms, this alignment process corresponds to adding subsequent indices in a single direction. As a result, the global configuration satisfies the relation $ x_1 < x_2 < \dots < x_d $ or $ x_1 > x_2 > \dots > x_d $. In either case, the resulting permutation remains optimal.

    In contrast, the second scenario highlights the potential consequences of incorrectly estimating the parameter $ a $. Let  $ i_j = 2j $. With this choice of $ J $, for every $ r > 1 $, $ \sum_{m=1}^{r-1} |D_{i_m} \cap D_{i_r}| = 1 $. Specifically, each subsequent neighborhood $ D_{i_r} $ is again disjoint from the neighborhoods $ D_{i_1}, \dots, D_{i_{r-2}} $ and shares exactly one neighbor (index $ i_r-1 $) with the neighborhood $ D_{i_{r-1}} $. At each step, the objective function \eqref{eq:rec_align} has two global minima: one at $ a = 1 $ and the other at $ a = -1 $. The parameter $ b $, which controls translation, is determined in the alignment process such that the locations of the shared neighbor ($ i_k-1 $) in the local configuration vectors $ x^{(k-1)} $ and $ x^{(k)} $ overlap.
    If at any step the procedure erroneously selects the value for $ a $, the direction of adding subsequent elements reverses. Until the next error occurs, elements continue to be added in the wrong direction. Consequently, due to the changes in direction caused by each additional error, the global configuration will exhibit regions affected by ``errors", where relatively distant elements are interspersed. Consequently, the width of $ \boldsymbol{\Sigma} $ under the resulting permutation $ \pi $ may be large.
\end{example}

The last criterion concerns the computational complexity of determining the parameters in the alignment process.
According to the estimation formula \eqref{eq:optimal_ab}, it is essential to identify which sets $ D_{i_m} $ intersect nontrivially with $ D_{i_r} $.
Since each set $ D_j $  intersects only with sets $ D_m $ for which $ m \in \bigcup_{i \in D_j} D_i $, we can construct the sequence $ i_r $ to leverage this relationship.
This approach reduces the number of intersections evaluated in \eqref{eq:optimal_ab}.

To present the main algorithm, we define the order-$t$ neighborhoods $D_i(t)$ and their outskirts $L_i(t)$ recursively. Initially, $D_i(0)=\{i\}$ and $L_i(0)=\{i\}$, and for $t\ge 1$
\begin{align*}
    L_i(t) & =\bigcup_{j\in L_i(t-1)}D_j\setminus D_i(t-1),\,\,
    D_i(t)=D_i(t-1)\cup L_i(t).
\end{align*}
Since $D_i(0) \subseteq D_i(1)\subseteq \dots \subseteq D_i(t)$ for any $t$, the outskirts $L_i(0) = D_i(0)$, $L_i(s) = D_i(s) \setminus D_i(s-1)$, for $s \ge 1$, are disjoint.
We refer to $L_i(s)$  as the outskirt of order $ s $ associated with row $ i $.
For $j \in L_i(t)$, we have
\begin{equation}
    \label{eq:subset}
    D_j \subseteq L_i(t-1) \cup L_i(t) \cup L_i(t+1).
\end{equation}
Consequently, if $j \in L_i(t)$, $m \in L_i(t+s) $ and $s\geq 3$, then $D_j \cap D_m = \emptyset$.

\begin{algorithm}
    \caption{Skeleton selection}
    \label{alg:skeleton}
    \begin{algorithmic}\vspace{2mm}

        \STATE{{\sc Input:} $D$ -- A list of $d$ subsets of $\{1,\dots, d\}$ representing neighborhoods} \vspace{2mm}
        \STATE{$ z = 1\,;\,t=0$} ;\vspace{1mm}
        \STATE{$J[z] = {\rm random} (1, \dots, d)$
        };  \COMMENT{{\tt \small Initialization of the skeleton, an integer at random}}\vspace{1mm}
        \STATE{$ L[[t]] =\{J[z]\} $};  \COMMENT{{\tt \small Initialization of the list of outskirts}}\vspace{1mm}
        \STATE{$ N[[t]] =\{J[z]\} $} ;  \COMMENT{{\tt \small Initialization of the list of neighborhoods}}\vspace{2mm}
        \WHILE{$L[[t]] \neq \emptyset$}
        \STATE{$t=t+1$} ;
        \STATE{$L[[t]]=\bigcup_{j\in L[[t-1]]}D[[j]]\setminus N[[t-1]]$
        }; \COMMENT{{\tt \small The t-th order outskirts of J[z]}}\vspace{1mm}\\
        \STATE{$N[[t]]=N[[t-1]]\cup L[[t]] $
        }; \COMMENT{{\tt \small The t-th order neighborhood of J[z]}}\vspace{1mm}\\
        \ENDWHILE \vspace{2mm}
        \FOR{$r=2$ to $t$}\vspace{1mm}
        \STATE{$X=L[[r]]$} ;\vspace{2mm}
        \WHILE{$X \neq \emptyset$}\vspace{1mm}
        \STATE{ $j = {\rm random} (X )$;\,$W  = D[[j]]\cap L[[r-1]]$}; \COMMENT{{\tt \small random pick from the outskirt}}\vspace{1mm}
        \STATE{$z = z+1$;\, $J[z] = {\rm random}\left(
            \underset{a \in W}{\operatorname{argmax}} |D_a|\right)$};\vspace{1mm}
        \STATE{ $X = X \setminus D[[J[z]]]$} ;\vspace{1mm}
        \ENDWHILE \vspace{1mm}
        \ENDFOR

        \vspace{2mm}
        \STATE{{\sc Return:} $J$ -- The skeleton of rows.  }
    \end{algorithmic}
\end{algorithm}

Algorithm~\ref{alg:skeleton} schematically presents the selection of $J$. During the selection process, we track to which $ L_i(t) $ each element belongs.
Using \eqref{eq:subset}, the number of sets in \eqref{eq:optimal_ab} can be efficiently reduced.
More specifically, for a fixed $ r $, if we partition the index set $ I = \{i_1, \dots, i_{r-1}\} $ into disjoint subsets $ I(t) = I \cap L_j(t) $ and assume that $ i_r \in L_j(s) $, then in formula~\eqref{eq:optimal_ab}, we only consider indices within the set $ I(s-2)\cup I(s-1)\cup I(s) \cup I(s+1) \cup I(s+2)$ instead of the entire set $ I $. For large $ r $, this can significantly reduce computational complexity.
Additionally, for any $u$ where $L_i(u) \neq \emptyset = L_i(u+1)$, the sequence of the sets  $ L_i(1), L_i(2), \ldots, L_i(u) $ is a partition of $\{1,\dots,d\}$. This guarantees the completeness of the skeleton. Reducing the size of $ J $ improves computational efficiency. While our method attempts to limit its size, our priority is ensuring accurate parameter estimation during the alignment process. Consequently, the algorithm is designed primarily with this objective in mind.

In Algorithm \ref{alg:skeleton}, randomness occurs in three places: in the selection of the initial index $ i_1 $, in the choice of an element $ j $ from the set $ X $, and in the selection of its neighbor $ i_z $. A key aspect in this context is to obtain a sequence $ i_r $ with the property that $ |D_{i_r}| \approx 2l^{(\pi)}_{i_r} +1 $ (for optimal $\pi$), as this yields a small upper bound for the approximation of distance matrices. Currently, there is no criterion for comparing candidate indices. Therefore, we select a random index with a neighborhood as large as possible in a set $W$ (see Algorithm \ref{alg:skeleton}), as this may contribute to a more accurate estimation of the parameter \( a \) (particularly when the neighborhood overlaps significantly with those of previously selected indices). Additionally, such a choice may facilitate a reduction in the number of indices in the set \( J \), since fewer indices might be required to achieve full coverage of the entire index set.

\section {Theoretical and simulation-based performance study}
\label{sec:sim}
The code required to reproduce the simulations and computations is publicly available in the GitHub repository \cite{packing_algorithm_github}.
We evaluate the performance of the packing algorithm on simulated permuted band and block-tridiagonal matrices, and to some extent, on dyadic matrices.
Note that band matrices are, in fact, an important subclass of block-tridiagonal matrices.
The half-bandwidth of a matrix $\boldsymbol \Sigma$ is given by $\lambda=\|\boldsymbol l\|_\infty = \max\{ l_i: i=1,\dots, d\}$, where $\boldsymbol l=(l_i)_{i=1}^d$ is defined as in Definition~\ref{def:band} under the identity permutation. Furthermore, we benchmark our packing algorithm against Variable Neighbourhood Search (VNS) \cite{mladenovicVariableNeighbourhoodSearch2010}, a well-established bandwidth-minimization heuristic that effectively searches for permutations that minimize the $\ell_\infty$-norm.
Since the original implementation of \cite{mladenovicVariableNeighbourhoodSearch2010} is not publicly available, we implement VNS following the algorithmic description in the paper. Our implementation of this algorithm can be found in \cite{packing_algorithm_github}.

Although our algorithm is primarily designed for permuted band and block-tridiagonal matrices, simulations indicate that it also has the potential to effectively detect dyadic structures.
This is because the algorithm tends to pack nonzero terms densely around the diagonal, thereby revealing group patterns in row-wise widths of the packed matrix.

We begin by examining the role of higher-order neighborhoods in the algorithm.
Let $\boldsymbol \Sigma$ be a $0$-$1$ matrix and $\boldsymbol \Sigma(t)$ denote the $0$-$1$ matrix representing row connections in $t$-th order neighborhoods, more precisely, $\boldsymbol \Sigma_{ij}(t) =1$ if and only if  $j \in D_i(t)$.
For a permuted full band matrix $\boldsymbol \Sigma$, we define $\tilde\Pi(\boldsymbol\Sigma)$ as the set of permutations $\pi$ such that $\boldsymbol\pi^\top \boldsymbol\Sigma \boldsymbol\pi$ is a full band matrix.
Recall that $\Pi(\boldsymbol \Sigma)$ is the set of optimal permutations, i.e., those that minimize the half-width.
Intuitively, for a permuted full-band matrix $\boldsymbol \Sigma$, one may expect $\tilde\Pi(\boldsymbol\Sigma) = \Pi(\boldsymbol\Sigma)$.
Similarly, since the matrices $\boldsymbol \Sigma(t)$ also form permuted full-band matrices, one may expect $\tilde\Pi(\boldsymbol\Sigma(t)) = \Pi(\boldsymbol\Sigma(t))$.
These results are formally established in the following two theorems, with proofs provided in the Appendix.

\begin{figure}
    \centering  \includegraphics[width=\linewidth,height=0.5\linewidth]{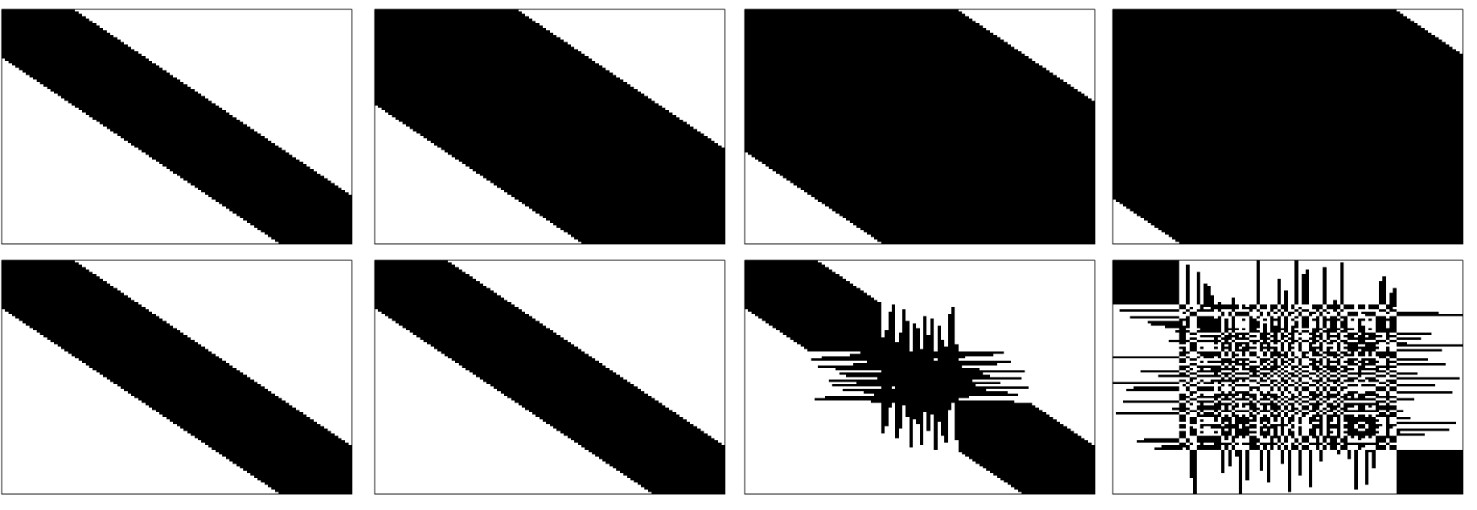}\vspace{-5mm}
    \caption{Effect of increasing the neighborhood size for finding the optimal permutation for the matrix $\boldsymbol{\Sigma}$. The first row shows   $\boldsymbol{\Sigma}(t)$, $t=1,2,3,4$ for a $100 \times 100$ full band matrix $\boldsymbol{\Sigma}$ with a half-bandwidth $\lambda = 20$. The second row shows the matrices recovered by applying the permutations obtained from the packing algorithm to the randomly permuted matrices in the first row.}
    \label{fig:grid}
\end{figure}

\begin{theorem}
    Let $\boldsymbol \Sigma $ be
    a $ d \times d $ permuted full band matrix. Then $\tilde \Pi(\boldsymbol{\Sigma}) \subseteq \Pi(\boldsymbol{\Sigma})$.
    \label{t:bandopt}
\end{theorem}

\begin{remark}
    A stronger result holds that, if $\lfloor \frac{d}{\lambda} \rfloor > \lfloor \frac{d-1}{\lambda+1} \rfloor$, then $\tilde \Pi(\boldsymbol{\Sigma}) = \Pi(\boldsymbol{\Sigma})$.
    The inequality $\lfloor \frac{d}{\lambda} \rfloor > \lfloor \frac{d-1}{\lambda+1} \rfloor$ is always satisfied when $ \lambda^2 \le d $. On the other hand, when $ \lambda^2 > d $, it may or may not hold. Simulations suggest that the inequality is primarily a technical condition and that, regardless of whether it holds, the equality $\tilde \Pi(\boldsymbol{\Sigma}) = \Pi(\boldsymbol{\Sigma})$ remains valid.
    However, when $\lfloor \frac{d}{\lambda} \rfloor = \lfloor \frac{d-1}{\lambda+1} \rfloor$, the result is still a conjecture.
\end{remark}
\begin{theorem}
    \label{T: per full band matr}
    Let $t\geq1$ satisfy $t\lambda\leq d$. If a $d\times d$ matrix $\boldsymbol{\Sigma}$ is such that for a permutation $\tilde \pi$, $\tilde{\boldsymbol \pi}^{\top} \boldsymbol \Sigma\tilde{\boldsymbol \pi}$ is a full band matrix with the half-bandwidth $\lambda \geq 1$, then the following conditions hold:
    \begin{itemize}
        \item[\it i)]  $\tilde{\boldsymbol \pi}^{\top}\boldsymbol \Sigma (t)\tilde{\boldsymbol \pi}$ is a  full band matrix with a half-bandwidth $t\lambda$.
        \item[\it ii)] If  $2\lambda t+1 \le d$, then $\tilde \Pi(\boldsymbol{\Sigma}(t))=\{\tilde \pi,  \rho\circ \tilde \pi\}$, where $\rho(i) = d+1-i$ for $i = 1,...,d$,
        \item[\it iii)] If  $2\lambda t + 1 >d$, then
              $|\tilde \Pi(\boldsymbol{\Sigma}(t))|=2 \cdot (2\lambda t + 2 - d)!$ and $\tilde \Pi(\boldsymbol{\Sigma}(t))\subset \tilde \Pi(\boldsymbol{\Sigma}(t+1))$.
    \end{itemize}

\end{theorem}

We highlight two key observations that follow from the above results.
First, the set of optimal permutations can be explicitly characterized, since any permutation $\pi$ leading to a full band $\boldsymbol{\pi}^\top \boldsymbol{\Sigma\pi}$ is explicitly constructible (see the proof of Theorem~\ref{T: per full band matr}). Moreover, simulations show that the packing algorithm applied to permuted full band matrices (for various $d$, $\lambda$) consistently recovers optimal permutations.
Second, {\it ii)} and {\it iii)} imply that the problems of finding optimal permutations for $\boldsymbol\Sigma$ and $\boldsymbol\Sigma(t)$ are equivalent when $2\lambda t + 1 \le d$—that is, $\tilde \Pi(\boldsymbol{\Sigma}) = \tilde \Pi(\boldsymbol{\Sigma}(t)) = \{\tilde \pi, \rho \circ \tilde \pi\}$. In this regime, both matrices can be used interchangeably, as confirmed by simulations (Figure~\ref{fig:grid}, first two columns).
When $2\lambda t + 1 > d$ (last two columns), this equivalence breaks down and the algorithm applied to $\boldsymbol\Sigma(t)$ need not recover the optimal permutation for $\boldsymbol\Sigma$.
These observations offer useful intuition: using $\boldsymbol{\Sigma}(t)$ can improve the approximation of the optimal permutation, but if $t$ is too large, this advantage may disappear, as seen in the third and fourth columns of Figure~\ref{fig:grid}.

\begin{figure}
    \centering
    \includegraphics[width=0.49\linewidth]{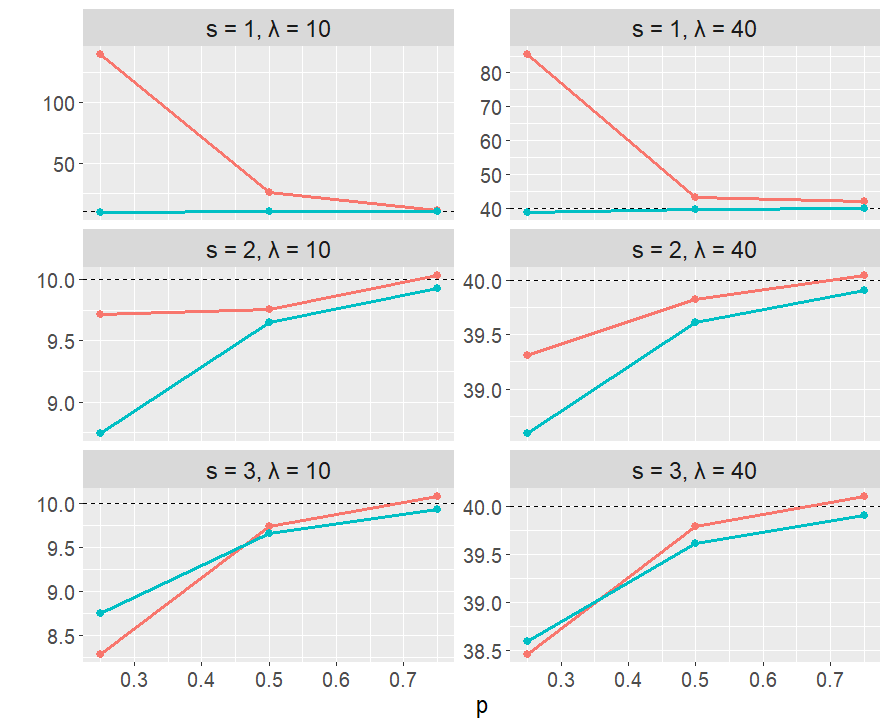}
    \includegraphics[width=0.49\linewidth]{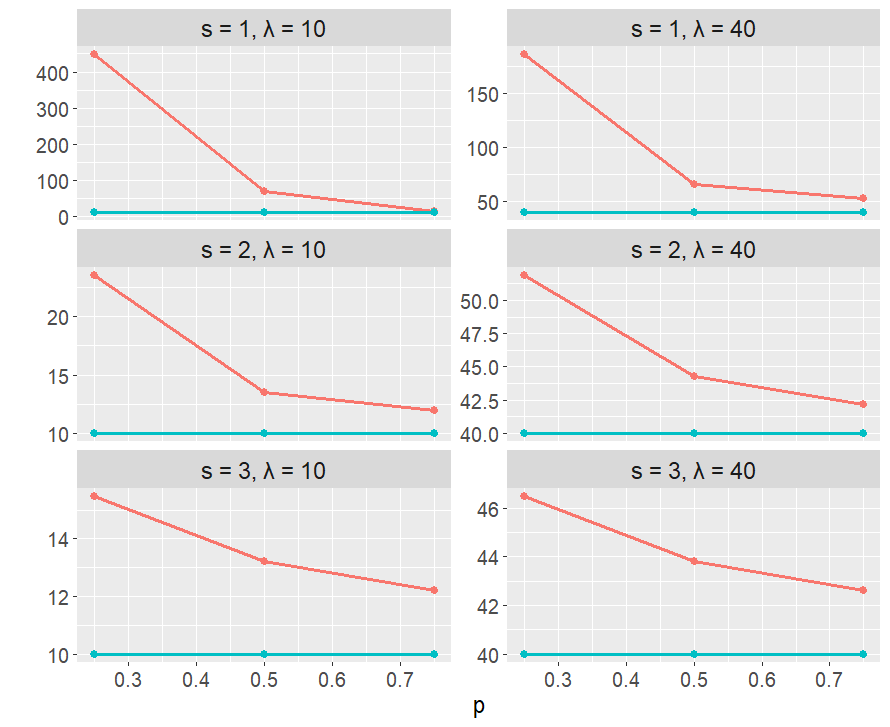}
    \vspace{-2mm}
    \caption{The first two columns show the average half-width \(\|\boldsymbol l^{(\pi)}\|_1/d\) for \(\lambda=10\) and \(40\), respectively; the last two show \(\|\boldsymbol l^{(\pi)}\|_\infty\) for the same values of \(\lambda\). Rows correspond to \(s\in\{1,2,3\}\). The blue lines in each plot represent values calculated over the initial (before permutation) matrices, while the red lines correspond to averages calculated over the final matrices obtained by the packing algorithm.}
    \label{fig:bandmatresults}
\end{figure}
\vspace{1.5mm}
\subsection{Permuted band matrices}
As shown in Theorems~\ref{t:bandopt} and~\ref{T: per full band matr}, the full band matrices form a restrictive but analytically tractable case.
Although theoretical guarantees are currently not known for more general settings, simulations suggest that the packing algorithm remains effective for incomplete but not too sparse band matrices.
This effectiveness is due to two key factors. First, the quality of approximation depends on the local density of $\pi(D_i)$ within $\pi(i) + C(l_i^{(\pi)})$. According to Remark~\ref{rem:imp}, the denser these sets, the tighter the bounds in Theorem~\ref{T:1.3}, and thus the better the approximation of $\mathbf{G}_\pi$ by $\mathbf{A}_{\boldsymbol{\Sigma}}$.
Second, the sizes of $D_i$'s affect the performance of CMDS.
A smaller $D_i$ leads to a smaller local distance matrix and may reduce accuracy and robustness.

Using higher-order neighborhoods often yields more accurate approximations because $\boldsymbol{\Sigma}(t)$ helps mitigate both factors that degrade the estimation quality.
First, increasing the neighborhood order expands the set of indices for which pairwise distances can be computed. This allows the restrictive condition $D_i \cap D_j \neq \emptyset$ (see the discussion of Theorem~\ref{T:1.3} assumptions) to be relaxed to $D_i(t) \cap D_j(t) \neq \emptyset$. This is particularly beneficial since the sets $D_i(t)$ are often significantly larger than $D_i$.
Second, for an optimal $\pi$, the matrix $\boldsymbol{\pi}^\top \boldsymbol{\Sigma}(t) \boldsymbol{\pi}$ typically has a higher proportion of nonzero entries within the band compared to $\boldsymbol{\pi}^\top \boldsymbol{\Sigma} \boldsymbol{\pi}$, and the efficiency of the algorithm observed for permuted full-band matrices kicks in.

In the simulations, we generated $1000 \times 1000$ symmetric band matrices, later referred to as initial matrices. Each was constructed by starting from a full band matrix with half-bandwidth $\lambda \in \{10, 40\}$ and independently setting entries above the diagonal within the band to 0 with probability $1 - p$, for $p \in \{0.25, 0.5, 0.75\}$ (diagonal elements were not zeroed). Symmetry was enforced by mirroring the upper-triangular part onto the lower triangle.
Each initial matrix was randomly permuted. For the resulting matrix, we computed $\boldsymbol{\Sigma}(s)$ with $s \in \{1, 2, 3\}$ and applied the packing algorithm to each of them to obtain a packed version of the matrix.
Performance was evaluated by comparing the half-width and half-bandwidth of final and initial matrices. For each combination of $(\lambda, p, s)$, the experiment was repeated 100 times.

The first two columns of Figure \ref{fig:bandmatresults} show simulation results for the half-width as a function of $p$, with blue and red curves representing the initial and final matrices, respectively. The blue curves align with intuition: the initial matrices are constructed by randomly zeroing elements within the band and the increased proportion of zeros results in a smaller half-width. This pattern is similar for both values of \(\lambda\), suggesting that $\lambda$ has a limited effect on the shape of the blue curves.
Each row corresponds to a different neighborhood order $s$, but since the initial matrices depend only on $p$ and $\lambda$, the blue curves remain invariant across rows. Thus, they serve as a reference for assessing the red curves, which represent the outcome of the packing algorithm applied to $\boldsymbol{\Sigma}(s)$ for various values of $s$.

We now examine the red curves representing results for final matrices. For $s = 1$, they exhibit a decreasing trend: accuracy improves with increasing $p$. In sparse cases ($p = 0.25$), deviations from the optimal are notable, but denser matrices yield more accurate solutions. This supports earlier findings, indicating that when sets $\pi(D_i)$ within $ \pi(i) + C(l^{(\pi)}_i) $ become denser, the approximation of the true distances between rows becomes more accurate, which enhances the quality of the estimated permutation.
As shown in the second and third rows, higher-order neighborhoods further improve the results, with red and blue curves closely aligned. Notably, in the third row for $p = 0.25$, the packing algorithm (on average) is able to find permutations resulting in a smaller average half-width than those of the initial matrices. Although the initial matrix may appear at first glance optimal, the removal of 75\% of nonzeros (for $p = 0.25$) creates a structure that paradoxically enables denser diagonal packing—highlighting the algorithm’s effectiveness in minimizing half-width.

\begin{figure}
    \centering
    \includegraphics[width=0.27\linewidth]{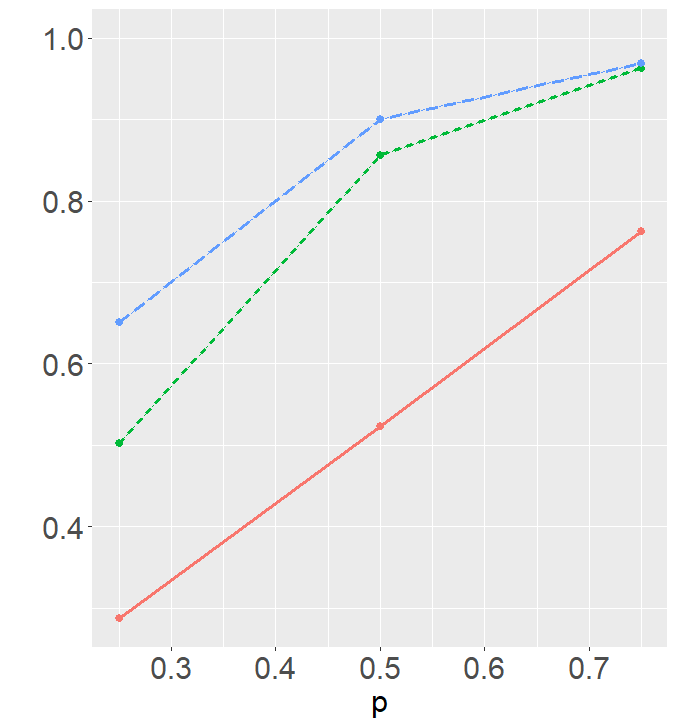}
    \includegraphics[width=0.6\linewidth]{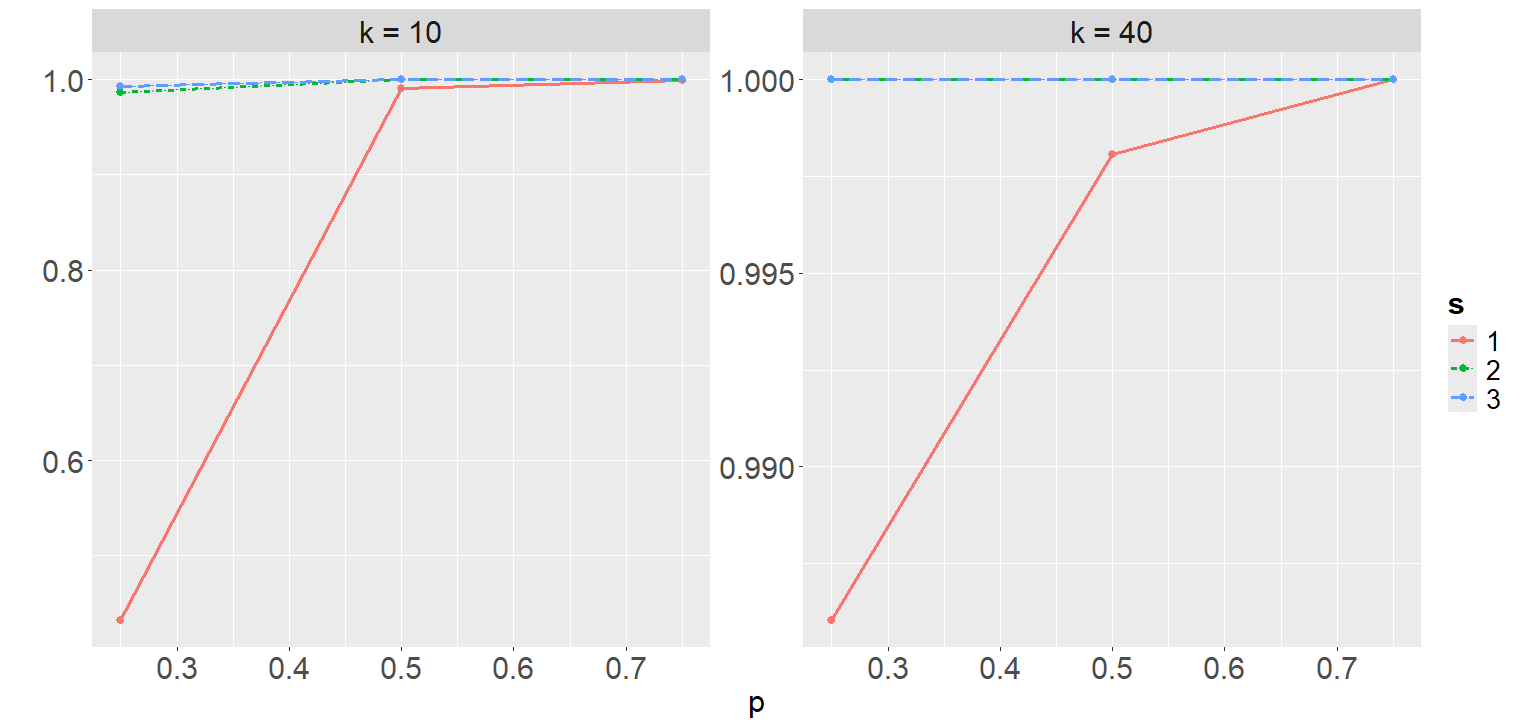}
    \vspace{-1mm}
    \caption{
        {\it Left:} the average $F(s)$ obtained for band matrices, $\lambda =10$.
            {\it Middle-Right:} the average proportion of nonzero elements within the block-tridiagonal structure returned by the packing algorithm applied to randomly permuted block-tridiagonal matrices with widths $10$ and $40$. The graphs are functions of the level of filling $p$ for  $s \in \{1,2,3\}$.
    }
    \label{fig:bandmatresults2}
\end{figure}

We note a close relation to the bandwidth minimization problem, which seeks a permutation $ \pi $ that minimizes the half-bandwidth $\Vert \boldsymbol{l}^{(\pi)} \Vert_{\infty}$ of $ \boldsymbol{\pi}^\top \boldsymbol{\Sigma} \boldsymbol{\pi} $.
Interestingly, although designed for a different purpose, the packing algorithm often produces near-optimal solutions for this problem, as shown in the third and fourth columns of Figure~\ref{fig:bandmatresults}.
As expected, the blue curves are constant at $\lambda$. While removing nonzeros could reduce the half-bandwidth, such cases are rare under the considered settings.
Higher-order neighborhoods consistently lead to results approaching the optimum.

$\boldsymbol{\Sigma}(s)$ must be a band matrix with a half-bandwidth no greater than $\min\{s\lambda, d\}$.
To compare band density across parameter settings, we define {\em the level of (band) filling}
$$
    F(s) = \frac{
        \sum_{i = 1}^d
        |\{j:\boldsymbol{\Sigma}_{ij}(s)\neq 0\}|
    }{
        \sum_{i = 1}^d|\{j:|i-j|\le s\lambda\}|
    },
$$
which quantifies the fraction of nonzero entries within the band. Higher filling yields a structure closer to a full band matrix, improving distance approximation and enhancing the quality of the resulting permutation.

\begin{figure}
    \centering
    \includegraphics[width=0.7\linewidth]{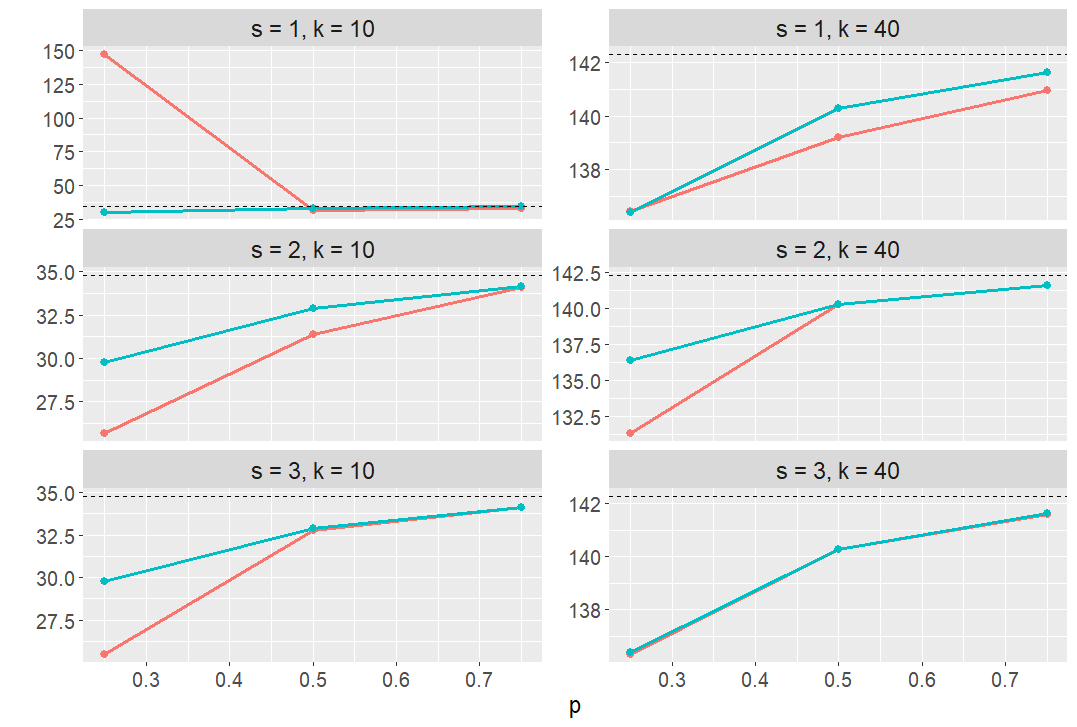}
    \vspace{-2mm}
    \caption{
        The average half-width $\Vert \boldsymbol l^{(\pi)} \Vert_1/d$ for permuted block-tridiagonal matrices is shown as a function of $p$ with $N =4$. Rows correspond to \(s\in\{1,2,3\}\), and columns correspond to \(k\in\{10,40\}\).
        The blue lines represent values calculated over the initial (not permuted) matrices, while the red lines correspond to averages calculated over the final matrices obtained by the packing algorithm.
    }
    \label{fig:tridblockmatresults3_1}
\end{figure}

Figure~\ref{fig:bandmatresults2} (left panel) presents the level of band filling as a function of $p$.
The solid line for $s = 1$ is linear, reflecting initial matrix construction. Its values slightly exceed $p$ due to the exclusion of diagonal elements from random zeroing. In contrast, the dashed lines ($s \in \{2, 3\}$) show noticeably higher filling levels, which enhance the algorithm’s performance by enabling better approximations of the optimal permutation.
We observe that the filling level increases with $s$, especially for $p = 0.25$. However, beyond a certain threshold, further increases in $F$ yield little improvement in the permutation quality. This is evident in the second and third rows of Figure \ref{fig:bandmatresults} for $p \geq 0.5$, where red and blue curves nearly coincide.
This observation is critical because increasing the neighborhood order $s$ increases the algorithm’s computational cost and, as stated in Theorem \ref{T: per full band matr}, when $2s\lambda + 1 > d$, finding the optimal permutation for $\boldsymbol{\Sigma}$ and for $\boldsymbol{\Sigma}(s)$ are no longer equivalent ($\tilde \Pi(\boldsymbol{\Sigma}) \subsetneq \tilde \Pi(\boldsymbol{\Sigma}(s))$). As a result, the algorithm may fail to return an accurate solution, as seen in Figure~\ref{fig:grid}.
\vspace{1.5mm}

\subsection{Permuted block-tridiagonal matrices}
\begin{figure}[!htbp]
    \centering
    \includegraphics[width=0.88\linewidth]{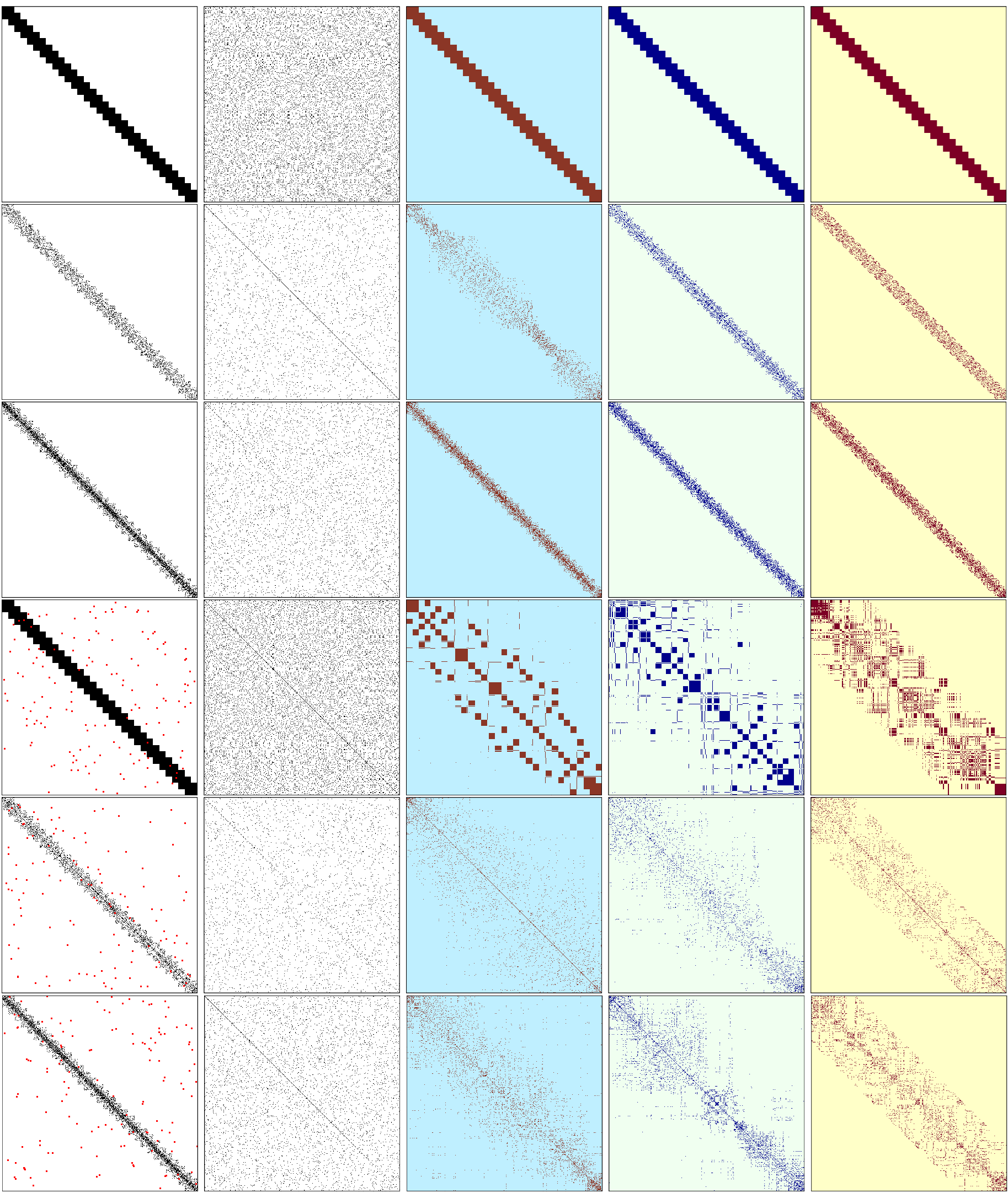}
    \vspace{-0mm}
    \caption{Packing methods for block-tridiagonal matrices and their noise-distorted versions.  {\it Column 1:} The matrices, from the top to the bottom: the fully-filled matrix, a sparse version, a denser sparse version near the diagonal, and the same three cases with noise (red dots).  {\it Column 2:} Randomly permuted matrices.  {\it Column 3:} The first-order neighbor reconstructions.  {\it Column 4:} The second-order neighbor reconstructions.  {\it Column 5:} The band matrices obtained by the VNS algorithm.
        \vspace{-3mm}
    }
    \label{fig:tridblock}
\end{figure}

Block-tridiagonal matrices constitute another subclass of dyadic matrices; see the top three plots in the first column of Figure~\ref{fig:tridblock}, where the fully and sparsely filled cases are presented.
We have also implemented a simulation setup analogous to that in the band-matrix study, except that the dimension $d = k(2^N - 1)$ is chosen to match complete dyadic matrices. We set $N = 4$ (height) and $k \in \{10, 40\}$ (block size).
Results are presented in Figure~\ref{fig:tridblockmatresults3_1}, organized analogously to Figure~\ref{fig:bandmatresults}, with the block size $k$ replacing the half-bandwidth $\lambda$.

For the denser matrices ($p \geq 0.5$), the results remain consistently close to optimal (the generating ordering) across all tested values of $s$ and $k$.
Interestingly, for $k = 40$, the outcomes based on $\boldsymbol{\Sigma}$ ($s = 1$) outperform those obtained using higher-order neighborhoods ($s \in \{2, 3\}$).
These findings highlight the strong performance of the packing algorithm and suggest that higher-order neighborhoods are not always necessary for a high-quality approximation of the optimal permutation when the within-structure nonzeros are dense enough. In the sparsest case ($p = 0.25$), for $k = 40$, the best packing was obtained for $s = 2$, though the permutations resulting from $s \in \{1, 3\}$ also provide good approximations. For the parameter setting $p = 0.25$, $s = 1$, and $k = 10$, the algorithm yields weaker results. Nevertheless, incorporating higher-order neighborhoods ($s \in \{2, 3\}$) significantly improves the quality of the approximation. In particular, the algorithm successfully identifies a permutation that concentrates nonzero elements closer to the diagonal more effectively than simply using the initial matrix.

For the dyadic algorithm, it is important that the packing algorithm preserves the block-tridiagonal structure, which can be measured by the average proportion of nonzero elements retained within the block-tridiagonal structure of the original full block-tridiagonal matrix.
This characteristic is shown in Figure~\ref{fig:bandmatresults2} (middle and right panels) on a scale from 0 to 1.
When higher-order neighborhoods are utilized ($s \in \{2, 3\}$), nearly all nonzero entries remain within the target structure. Only for $p=0.25$ does a small fraction (below 3\%) fall outside. For $s = 1$, the structure is well preserved when $p \geq 0.5$, but results degrade for sparser matrices ($p = 0.25$).

\vspace{1.5mm}

\subsection{Robustness study}
Using the block-tridiagonal matrices and their randomly perturbed versions shown in Figure~\ref{fig:tridblock}, we compare the robustness of our packing algorithms (plots in the third and fourth columns) with that of the VNS algorithm (plots in the last column).
Three cases are considered: fully filled, sparsely and uniformly filled, and unevenly filled with greater density near the diagonal.
For the fully filled case, both algorithms recover the block-tridiagonal structure perfectly, as seen in the top row.
The reconstructions for the sparsely filled case are also very accurate, with some caveats.
The reconstruction based on first-order neighborhoods ({\it Column 3}) is slightly less compact. Incorporating second-order neighborhoods ({\it Column 4}) yields an almost exact reconstruction. Packing into a band matrix (VNS, {\it Column 5}) somewhat distorts the structure by ``cutting corners'' that are far from the diagonal.
For a sparse block-tridiagonal matrix with a denser distribution close to the diagonal, all reconstructions are rather accurate, with our algorithm retrieving the original structure more faithfully.

\begin{figure}[H]
    \centering
    \includegraphics[width=0.75\linewidth]{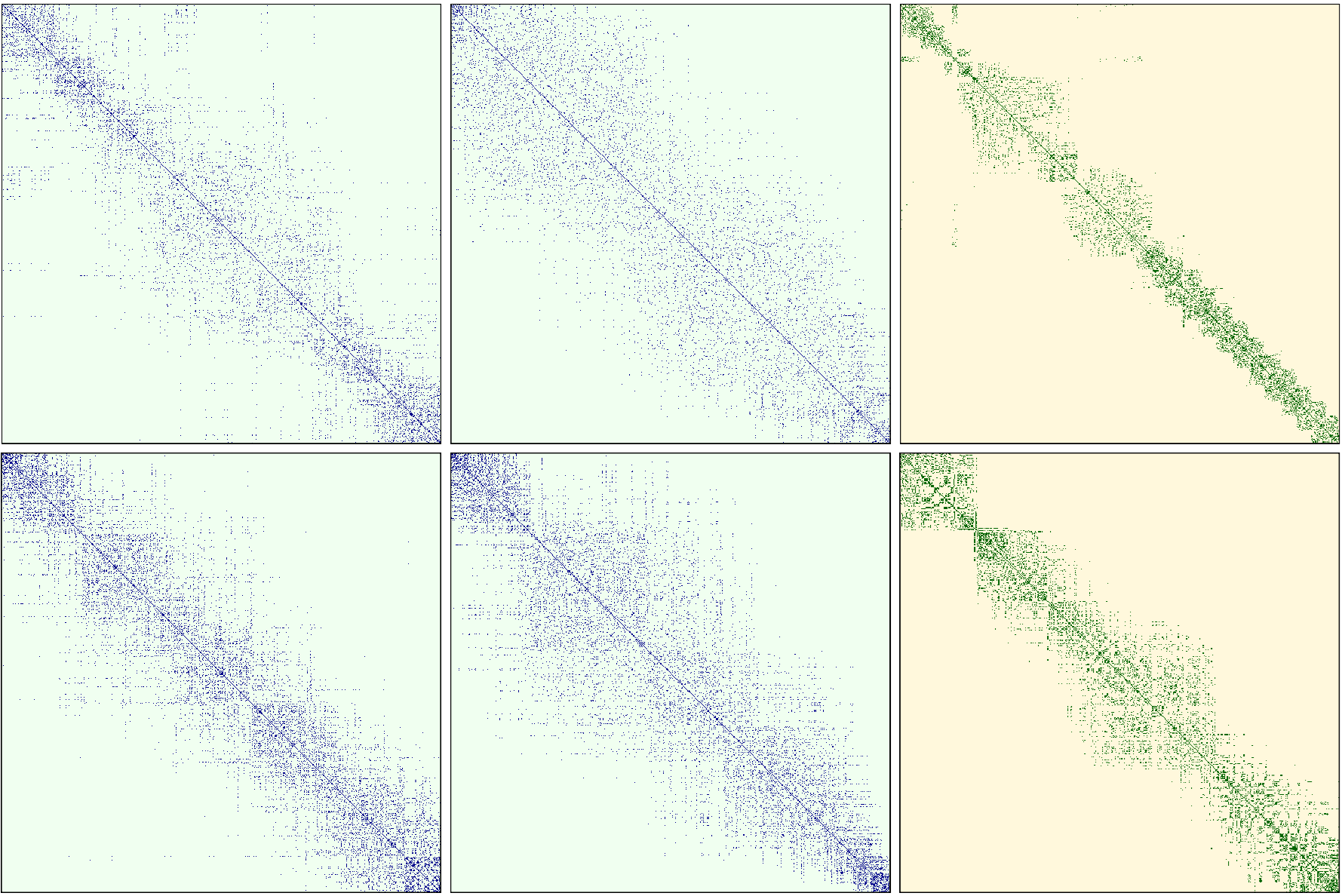}
    \caption{
        Reconstructions using third-order (left) and fourth-order (middle) neighborhoods for the last two cases in Figure~\ref{fig:tridblock} (the fifth and sixth rows). The right-hand column shows the reconstructions obtained by recursively removing approximately 100 outlying entries. Second-order neighborhoods are used in the top reconstruction and third-order neighborhoods in the bottom reconstruction.
    }
    \label{fig:tridblock_iterative}
\end{figure}

Our algorithm, which is based on the $\ell_1$-norm, can recover the structure of randomly perturbed matrices, a feature that may not be available to the bandwidth minimization method relying on the $\ell_\infty$-norm.
Due to the noise, it is more difficult to pack matrices tightly around the diagonal. However, for our algorithm, we observe that group structures emerge, especially when the second-order neighborhoods are used {\it (Column 4)}, while for the VNS algorithm, the original groupings are sliced and distributed arbitrarily within the band.
These structures become even more identifiable if the third- and fourth-order neighborhoods are considered. See the left-hand and middle plots in Figure~\ref{fig:tridblock_iterative}.
The final $\ell_1$-norms for the right-hand matrices are 25,915 and 39,743, respectively, compared with 93,839 and 86,335 for the third-order-neighborhood reconstructions.

In addition to a better reconstruction of the structural blocks, we observe identifiable isolated points around the groups.
These can be identified as outlying entries if the distance to the diagonal is used as a metric, and they likely impose challenges on a tight packing.
This feature is exploited in a simple iterative method based on the removal of outliers farthest from the diagonal. At each step, a small portion of the outliers is removed, and the algorithm is run on the remaining elements. The recursion is stopped when further reductions in the \(\ell_1\) criterion become small; in the simulated examples, this occurs after approximately 100 outlying entries have been removed.
The outcome of the implemented iterative procedure is presented in the right-hand panels of Figure~\ref{fig:tridblock_iterative}, which show a significant improvement in structural reconstruction.
It is worth pointing out that no analogous iterative VNS implementation is currently feasible.

\subsection{Permuted dyadic matrices}
\begin{figure}[!htbp]
    \centering
    \includegraphics[width=0.45\linewidth]{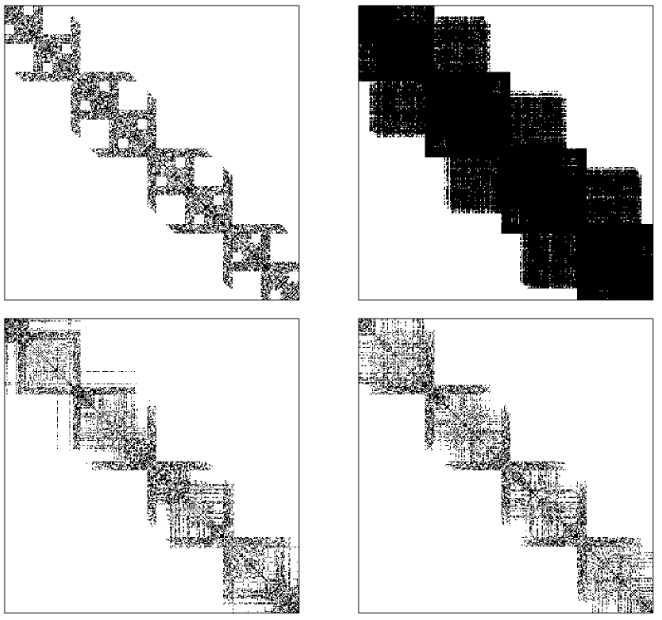}
    \hspace{2mm}\includegraphics[width=0.45\linewidth]{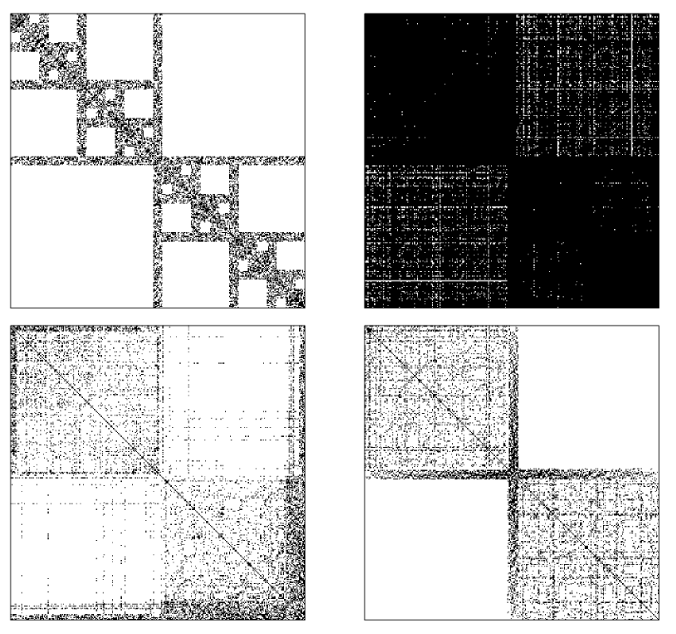}
    \caption{\it The performance of the packing algorithm for a banded dyadic matrix with $\lambda=60$ (the first two columns) and a dyadic matrix (the two remaining columns), with $N = 5, k = 10, p = 0.5$. The top row displays the original (non-permuted) matrices along with their second-order neighborhood matrices. The bottom row shows the matrices recovered by the algorithm.
    }
    \label{fig:alg bdm}
\end{figure}

Dyadic matrices have a recursive block-diagonal structure with sparsely populated nonzero blocks. As a result, our packing algorithm can be applied iteratively. To illustrate this approach, consider a banded dyadic matrix obtained via element-wise multiplication ($\mathbf{A}_{ij} = \mathbf{B}_{ij}\mathbf{C}_{ij}$ for every pair $ i,j $) of a full band matrix $\boldsymbol{B}$ and a dyadic matrix $\boldsymbol{C}$. See the upper-left panel of Figure~\ref{fig:alg bdm}, in which the first two columns show the application of the packing algorithm to a banded dyadic matrix, while the third and fourth correspond to a dyadic matrix.


For the banded dyadic example, packing both $\boldsymbol\Sigma$ and $\boldsymbol\Sigma(2)$ recovers the outermost block structure, with the second-order neighborhood matrix giving the clearer result because it more closely resembles a full block-tridiagonal matrix. The internal dyadic structure is not recovered at this stage because it has little influence on the global half-width criterion. The algorithm's success with banded dyadic matrices is linked to their structural similarity in sparsity pattern to block-tridiagonal matrices. However, the algorithm remains competitive when the original matrix deviates from this structure. For the general dyadic matrix, direct packing produces a large half-width but reveals groups of strongly connected indices. Packing $\boldsymbol\Sigma(2)$ recovers the outermost block structure because blocks nearer the diagonal are more densely populated than those farther away.

\begin{figure}[!htbp]
    \centering
    \includegraphics[width=0.7\textwidth]{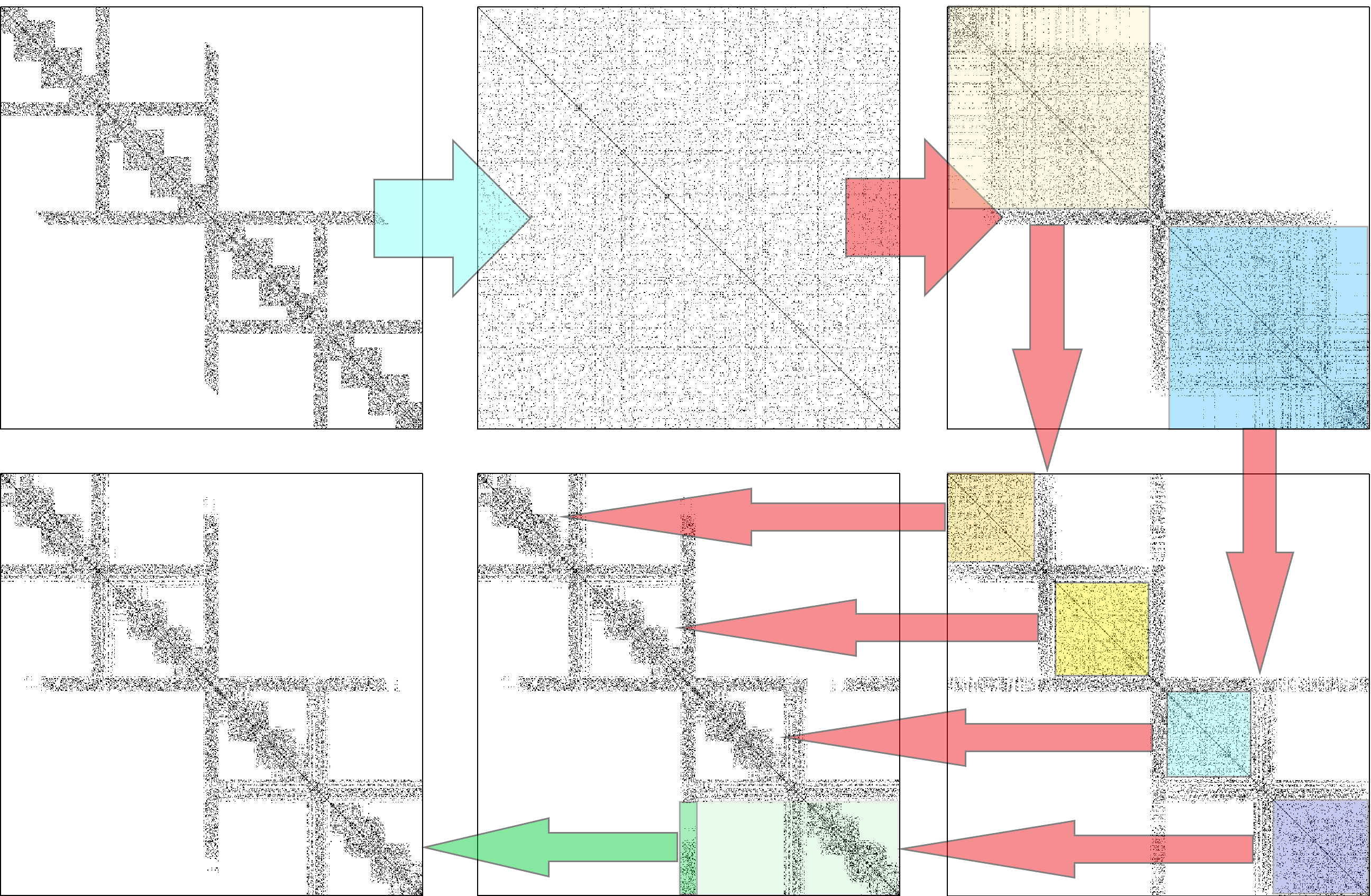}
    \caption{Example illustrating iterative application of the packing algorithm for determining the dyadic structure.}
    \label{fig:iter_pack}
\end{figure}

The algorithm’s ability to detect the outermost structure of a banded or non-banded dyadic matrix enables an iterative scheme for approximating the full reconstruction of the original dyadic form.
Figure \ref{fig:iter_pack} schematically illustrates how the packing algorithm can be applied iteratively in order to reconstruct a matrix with a dyadic structure.
In the top-left panel, a matrix that incorporates all matrix types considered so far is displayed. It was generated using the following parameters:
$k = 20$, $N = 5$, $p = 0.25$, $\lambda = 250$.
Subsequently, selected regions of the matrix were zeroed out in order to obtain several distinct structural patterns.
The entire matrix is of the \textbf{banded dyadic matrix} type. This is easily observed from the truncated arms of the main cross located at the center of the matrix. This main cross separates the matrix into two submatrices, the upper-left and the lower-right, which are disconnected from each other except through the main cross. Each of these submatrices is a \textbf{(non-banded) dyadic matrix} in the sense that every element on its main cross is equal to 1 with probability $p = 0.25$ (without truncation). Analogously, the main crosses inside these submatrices further split them into two additional submatrices, and each of the four resulting submatrices is of \textbf{block-tridiagonal} type.

After random permutation, applying the packing algorithm to $\boldsymbol\Sigma(2)$ recovers the main cross but not the internal ordering of the two diagonal submatrices. We therefore partition the rows into the main-cross indices and the two diagonal groups and apply the packing algorithm separately to the latter using second-order neighborhoods. Repeating this procedure on the resulting four submatrices recovers their internal structures. One submatrix is returned in reverse order; since $\pi$ and $\rho\circ\pi$ are equivalent for the dyadic algorithm, the final reversal is applied only for visual consistency.

It remains to describe the criterion used for partitioning the rows into three subsets. In the example, we employ an intuitive, simple rule which consists of identifying the largest possible rectangle filled entirely with zeros in the top-right corner of the reordered matrix.
Let \(\mathbf A\) denote the reordered matrix at the current stage. We identify the largest all-zero rectangle in its upper-right corner by selecting a pair \((i,j)\) such that
\[
    A_{i'j'}=0
    \quad\text{for all }i'\leq i\text{ and }j'\geq j,
\]
and maximizing \(f(i,j)=i(d+1-j)\) over all such pairs.
Once this pair is selected, the partitioning proceeds as follows:
\begin{itemize}
    \item the rows \(\{1,2,\ldots,i\}\) are associated with the upper-left square,
    \item the rows \(\{i+1,i+2,\ldots,j-1\}\) are associated with the main cross,
    \item the rows \(\{j,j+1,\ldots,d\}\) are associated with the lower-right square.
\end{itemize}
The above criterion is relatively simple but not optimal. Both the partitioning criterion and the iterative packing procedure remain subjects of ongoing investigation and future development.


\section{Conclusions}
\label{sec:conclusion}
We propose an efficient approach for the factorization and inversion of sparse dyadic matrices, providing a mathematical formalization of the classical nested dissection method. For a fixed block size, the dyadic factorization has quasi-linear computational complexity in the matrix dimension. In addition, we introduce a new packing algorithm based on the $\ell_1$-norm and investigate its properties both theoretically and numerically. The algorithm is primarily designed for band and block-tridiagonal matrices, which form two important subclasses of dyadic matrices.

Owing to its effectiveness for block-tridiagonal matrices, the packing algorithm can also be applied iteratively to reveal and reconstruct the outermost structure of a general dyadic matrix. Several examples of such iterative procedures were presented and exhibited promising behavior, although a more detailed investigation of these methods remains a topic for future work.

While bandwidth minimization is not the primary focus of this work, our numerical experiments suggest that the proposed packing algorithm constitutes a novel and promising heuristic for this task. In particular, many efficient bandwidth minimization methods rely on an initial ordering that is already reasonably close to optimal; in this setting, our proposed algorithm may serve as an effective warm-start strategy for subsequent refinement procedures.
The present work emphasizes conceptual unification rather than exhaustive algorithmic optimization. Extending the framework to specialized graph classes and incorporating state-of-the-art heuristics remains an interesting direction for future research.


\section*{Disclosure statement}
The authors report there are no competing interests to declare.

\section*{Funding}
The second and third authors have been partially supported by the Swedish Research Council (VR) Grant DNR: 2020-05168.

\section*{Declaration of generative AI use}
The authors have used ChatGPT 5.6 Sol and Gemini 3.1 Pro for reference management, idea exploration, copyediting and language refinement.







\bibliographystyle{tfnlm}
\bibliography{references}

\appendix
\section{Additional Technical Results}
\label{sec:additional-theory}

\begin{proof}[Justification of Algorithm~\ref{alg:dyadic}]
    Let $\mathbf P_1= \mathbf P$, where $\mathbf P$ is obtained in the first step of the algorithm before the for-loop and $\mathbf P_l$ are $\mathbf P$ matrices obtained within the for loop, $l=2,\dots, N$.
    We argue that $\mathbf P_l$ orthonormalizes $\tilde{\mathbf  F}_l=\begin{bmatrix}
            \mathbf F_1 & \dots & \mathbf  F_l
        \end{bmatrix}$. The proof proceeds by induction. It is straightforward to verify that $\mathbf P_1$ orthonormalizes $\mathbf F_1$ by the definition of $\mathbf G_1(\cdot)$.
    Assume that for $l\ge 2$, $\mathbf P_{l-1}$ orthonormalizes $\tilde{\mathbf F}_{l-1}$.
    Observe that
    \begin{equation*}
        \tilde{ \mathbf F}_l{\mathbf P}_l= \begin{bmatrix}
            \tilde{\mathbf F}_{l-1} & \mathbf  F_{l}
        \end{bmatrix}
        \begin{bmatrix}\mathbf P_{l-1} & -\mathbf A \mathbf G \\
               \mathbf 0       & \mathbf G
        \end{bmatrix}=
        \begin{bmatrix}\tilde{\mathbf F}_{l-1}\mathbf P_{l-1} & \left(\mathbf F_l- \tilde{\mathbf F}_{l-1}\mathbf A
               \right)  \mathbf G
        \end{bmatrix}.
    \end{equation*}
    By induction,
    $
        \mathbf P_{l-1}^\top
        \tilde{\mathbf F}_{l-1}^\top
        \tilde{\mathbf F}_{l-1}
        \mathbf P_{l-1}=\mathbf I
    $,
    so it suffices to show
    $\tilde{\mathbf F}_{l-1}\mathbf P_{l-1} \bot \left(\mathbf F_l -\tilde{\mathbf F}_{l-1}\mathbf A
        \right)$ and that $\mathbf G$ orthonormalizes $\left(
        \mathbf F_l-\tilde{\mathbf F}_{l-1}\mathbf A  \right)$.
    The first part follows from
    \begin{equation*}
        \mathbf P_{l-1}^\top
        \tilde{\mathbf F}_{l-1}^\top
        \left(\mathbf F_l -\tilde{\mathbf F}_{l-1}\mathbf A\right)
    \end{equation*}
    \begin{equation*}
        = \mathbf P_{l-1}^\top
        \tilde{\mathbf F}_{l-1}^\top\mathbf F_l
        -
        \mathbf P_{l-1}^\top
        \tilde{\mathbf F}_{l-1}^\top
        \tilde{\mathbf F}_{l-1}
        \mathbf P_{l-1}
        \mathbf P_{l-1}^\top
        \begin{bmatrix}
            \boldsymbol \Sigma_{1,l} \\
            \vdots                   \\
            \boldsymbol \Sigma_{l-1,l}
        \end{bmatrix}=
        \mathbf P_{l-1}^\top
        \begin{bmatrix}
            \boldsymbol \Sigma_{1,l} \\
            \vdots                   \\
            \boldsymbol \Sigma_{l-1,l}
        \end{bmatrix}
        -
        \mathbf P_{l-1}^\top
        \begin{bmatrix}
            \boldsymbol \Sigma_{1,l} \\
            \vdots                   \\
            \boldsymbol \Sigma_{l-1,l}
        \end{bmatrix}.
    \end{equation*}

    The second part follows from
    \begin{align*}
         & \hphantom{=\ }\left(
        \mathbf F_l-\tilde{\mathbf F}_{l-1}\mathbf A  \right)^\top\left(
        \mathbf F_l-\tilde{\mathbf F}_{l-1}\mathbf A  \right) \\
         & =
        \boldsymbol \Sigma_{l,l}-
        \mathbf A^\top\tilde{\mathbf F}_{l-1}^\top
        \mathbf F_l-
        \left(
        \mathbf F_l-\tilde{\mathbf F}_{l-1}\mathbf A  \right)^\top\tilde{\mathbf F}_{l-1}\mathbf P_{l-1}
        \mathbf P_{l-1}^\top
        \begin{bmatrix}
            \boldsymbol \Sigma_{1,l} \\
            \vdots                   \\
            \boldsymbol \Sigma_{l-1,l}
        \end{bmatrix}
        \\
         & =
        \boldsymbol \Sigma_{l,l}-
        \begin{bmatrix}
            \boldsymbol \Sigma_{1,l}^\top ,
            \cdots ,
            \boldsymbol \Sigma_{l-1,l}^\top
        \end{bmatrix} \mathbf P_{l-1}\mathbf P_{l-1}^\top
        \begin{bmatrix}
            \boldsymbol \Sigma_{1,l} \\
            \vdots                   \\
            \boldsymbol \Sigma_{l-1,l}
        \end{bmatrix}
        = \tilde{\boldsymbol \Sigma},
    \end{align*}
    and $\mathbf G$ orthonormalizes vectors whose Gramians equal to $\tilde{\boldsymbol \Sigma}$.
\end{proof}

\begin{proof}[Proof of Lemma~\ref{lem:IHDED}]
    Since $\mathbf H$ is of shape $k(2^N - 2^{N-l+1})\times k(2^N - 2^{N-l+1})$, and $\mathbf E$ is of shape $k(2^N - 2^{N-l+1})\times k2^{N-l}$, $\mathbf {HE}$ has the same shape as $\mathbf E$ and by (\ref{eq:decompE}) and (\ref{eq:decompVH})
    \begin{equation}
        \label{eq:HE1}
        \mathbf {HE} = \left(\sum_{s=1}^{2^{N-l+1}}\mathbf H_s\right)\left(\sum_{r=1}^{2^{N-l}} \mathbf E_r \right)
    \end{equation}

    As we know, $\mathbf H_s$ is $\mathcal{HD}(l-1)$ in the $s$-th diagonal $k2^{l-1}\times k2^{l-1}$ block and zero otherwise.
    On the other hand, matrix $\mathbf E_r$ has a block-sparsity pattern corresponding to $\mathcal{ED}(N, l)_r$. Consequently, $\mathbf H_s \mathbf E_r$ is nonzero only if $s = 2r-1$ or $s = 2r$. Therefore, \eqref{eq:HE1} simplifies to
    $      \mathbf{HE} = \sum_{r=1}^{2^{N-l}} \left(\mathbf H_{2r-1} + \mathbf H_{2r}\right)\mathbf E_r$.
    Note that the nonzero part of $\mathbf H_{2r-1} + \mathbf H_{2r}$ is a block-diagonal matrix of size $k(2^l-2)\times k(2^l-2)$ with its two diagonal blocks being $\mathcal{HD}(l-1)$. Meanwhile, the nonzero part of $\mathbf E_r$ is of size $k(2^l-2)\times k$ and can be viewed as a $2\times 1$ block matrix whose blocks are each of size $k(2^{l-1}-1)\times k$. Multiplying these nonzero parts yields a result of size $k(2^{l}-2)\times k$, which is exactly the size of the nonzero block in $\mathbf E_r$. Hence, the product $\left(\mathbf H_{2r-1} + \mathbf H_{2r}\right)\mathbf E_r$ retains the same block-sparsity pattern as $\mathbf E_r$, and when summed over $r$, the resulting matrix $\mathbf{HE}$ exhibits the same block-sparsity pattern as $\mathbf E$, namely $\mathcal{ED}(N, l)$.

    In terms of complexity, according to Proposition~\ref{prop:HA}, the flops required to calculate each summand is $O\left((l-1)2^{l-1}k^3\right)$, hence the total flops required to obtain $\mathbf{HE}$ is
    $
        O\left((l-1)2^{N-1}k^3\right) = O\left(ldk^2\right).
    $
    The proof for $\mathbf {VE}$ is analogous. $\mathbf{VE}$ has the same block-sparsity pattern as $\mathbf E$. Based on Proposition~\ref{prop:VA}, the complexity is also $O(ldk^2)$.

    For the second part, according to \eqref{eq:decompE},
    $
        \mathbf E^\top\mathbf E =\left(\sum_{s=1}^{2^{N-l}} \mathbf E_s^\top\right)\left(\sum_{r=1}^{2^{N-l}} \mathbf E_r\right)$.
    From the block-sparsity pattern \eqref{eq:ED}, we have that $\mathbf E_s^\top\mathbf E_r$ is nonzero only if $s = r$. Hence,
    $
        \mathbf E^\top\mathbf E = \sum_{s=1}^{2^{N-l}} \mathbf E_s^\top\mathbf E_s$.
    To obtain each of the summands one essentially calculates multiplication of a $k\times k\left(2^{l}-2\right)$ nonzero matrix with its $k\left(2^{l}-2\right)\times k$ nonzero counterpart. This results in a $k\times k$ nonzero matrix, and it costs $O\left(k^3 2^l\right)$ flops.
    This means that $\mathbf E_s^\top \mathbf E_s$ results in a $k 2^{N-l}\times k2^{N-l}$ block diagonal matrix with block size $k\times k$ and only the $s$-th diagonal block of it is nonzero. Hence, $\mathbf E^\top \mathbf E$ is a $k 2^{N-l}\times k2^{N-l}$ block diagonal matrix with each of its diagonal block being nonzero. In other words, $\mathbf E^\top \mathbf E\in \mathcal D(N, l)$. The total number of flops required to compute $\mathbf E^\top \mathbf E$ is
    $
        O\left(k^3 2^N\right) = O\left(dk^2\right).
    $

    Finally, we note that $
        \mathbf D = \sum_{r=1}^{2^{N-l}} \mathbf D_r$,
    where $\mathbf D_r$ is a block diagonal matrix with $k\times k$ diagonal blocks and only its $r$-th diagonal block is nonzero.
    Again, referring to \eqref{eq:bspE}, we know that $\mathbf E_s \mathbf D_r$ is nonzero only if $s = r$. As a result,
    $
        \mathbf{ED} =  \sum_{s=1}^{2^{N-l}} \mathbf E_s\mathbf D_s$.
    To obtain each of the summands is essentially to calculate multiplication of a $k\left(2^{l}-2\right)\times k$ nonzero matrix with its $k \times k$ nonzero counterpart, resulting in a $k\left(2^{l}-2\right) \times k$ nonzero matrix. This means that $\mathbf E_s\mathbf D_s$ will have the same block-sparsity pattern as $\mathbf E_s$, and hence $\mathbf {ED}$ will have the same block-sparsity pattern as $\mathbf E$, which is $\mathcal{ED}(N, l)$.

    In terms of computational complexity, computing each summand requires $O(k^3 2^{l})$ flops, and in total, the number of flops required by $\mathbf{ED}$ is $O\left(k^3 2^N\right) = O\left(dk^2\right)$.
\end{proof}

\begin{proof}[Proof of Lemma~\ref{lem:IHDED_band}]
    Since the refined patterns $\widetilde{\mathcal{ED}}(N, l)$ and $\widetilde{\mathcal{ED}}^\prime(N, l)$ are subsets of $\mathcal{ED}(N, l)$, it follows from Lemma~\ref{lem:IHDED} that $\mathbf{VHE}$ inherits the sparsity pattern $\mathcal{ED}(N, l)$.

    From the proof of Lemma~\ref{lem:IHDED}, $
        \mathbf{HE} = \sum_{r=1}^{2^{N-l}} \left(\mathbf H_{2r-1} + \mathbf H_{2r}\right)\mathbf E_r$,
    which is essentially the sum of the products of a $k(2^l-2)\times k(2^l-2)$ block-diagonal matrix (with its two nonzero diagonal blocks having structure $\mathcal{HD}(l-1)$) and a $k(2^l-2)\times k$ submatrix whose block-sparsity pattern assures that only the rows in the $\mathcal{HD}(l-1)$ blocks corresponding to the nonzero entries specified by the block-sparsity pattern $\widetilde{\mathcal{ED}}^\prime(N, l)$ are involved in the multiplication. Thus $\mathbf{HE} \in \widetilde{\mathcal{ED}}^\prime(N, l)$. In terms of complexity, the number of flops required to evaluate $\mathbf{HE}$ is
    $
        O\left((l-1)k^32^{N-l+1}\right) = O(k^32^N) = O(dk^2).
    $

    Similarly, the structure of $\mathbf V\in \mathcal{IVD}(N,l)$ ensures that, in the product $\mathbf{VHE}$, only the columns corresponding to the nonzero entries in $\widetilde{\mathcal{ED}}^\prime(N, l)$ are involved. The  cost of this multiplication is
    $
        O\left(k^32^{N-l+1} \sum_{r=1}^{l-1} (2^{l-r}-1)\right) = O(k^3 2^N) = O(dk^2).
    $
\end{proof}

The proof of Theorem~\ref{T:1.3} requires the following auxiliary results.

\begin{lemma} \label{lemma 1.2}
    Under the assumptions of Theorem~\ref{T:1.3}, it holds
    \begin{equation*}
        \left |\pi(i) -\pi(j) \right| = \frac{1}{2}\left | \left (\pi(i)+ C(l_i^{(\pi)})\right ) \triangle \left (\pi(j)+ C(l_j^{(\pi)}) \right ) \right | ,
    \end{equation*}
    where $C(t) = \{-t,\dots, -1, 0, 1, \dots ,t\}$.
\end{lemma}
\begin{proof}
    Without loss of generality, one can assume that $\pi(i)<\pi(j)$.
    The assumed inequalities
    yield
    \begin{equation}
        \label{eq:ineq_B1}
        \pi(i) - l^{(\pi)}_i  \le \pi(j) - l^{(\pi)}_j,  \,\,
        \pi(i) + l^{(\pi)}_i  \le \pi(j) + l^{(\pi)}_j,   \,
        \,
        \pi(i) + l^{(\pi)}_i  \geq \pi(j) - l^{(\pi)}_j
    \end{equation}
    The first inequality shows that the minimum of set $ \pi(i) + C(l_i^{(\pi)}) $ bounds set $ \pi(j) + C(l_j^{(\pi)}) $ from below.
    The second indicates that the maximum of set $ \pi(j) + C(l_j^{(\pi)}) $ bounds set $ \pi(i) + C(l_i^{(\pi)}) $ from above.
    The third inequality implies that the two sets have a non-empty intersection.
    Consequently,
    \begin{multline*} \frac 1 2\left | \left (\pi(i)+ C(l_i^{(\pi)})\right ) \triangle \left (\pi(j)+ C(l_j^{(\pi)})\right ) \right | =\\
        =\frac 1 2\left [(\pi(j) - l^{(\pi)}_j) - (\pi(i) - l^{(\pi)}_i)\right ] + \frac 1 2\left [(\pi(j) + l^{(\pi)}_j) - (\pi(i) + l^{(\pi)}_i)\right ]
        = \pi(j) - \pi(i).
    \end{multline*}
    The equality between the first and second quantity follows from the fact that $\pi(k)+ C(l_k^{(\pi)})$, for $k = i,j$, are sets of consecutive integer numbers and inequalities~\eqref{eq:ineq_B1}.
\end{proof}

\begin{lemma} \label{lemma 2.1}
    Suppose $Z_1,Z_2,A_1,A_2$ are subsets of $\{1,\dots,d\}$.
    For any permutation $\pi$ of $(1,\dots,d)$, if  $\pi(Z_1) \subseteq A_1$ and $\pi(Z_2) \subseteq A_2$, then
    \[
        \vert Z_1 \cap  Z_2 \vert \le \vert A_1 \cap A_2 \vert
        \le \vert Z_1 \cap  Z_2 \vert
        +\vert  A_1 \vert - \vert  Z_1 \vert + \vert A_2  \vert - \vert   Z_2 \vert.
    \]
\end{lemma}
\begin{proof}
    By the assumption, $\vert \pi(Z_1) \cap  \pi(Z_2) \vert \le \vert A_1 \cap A_2 \vert$.
    On the other hand
    $
        \pi(Z_1) \cap  \pi(Z_2) = \{\pi(m): m \in Z_1\} \cap \{\pi(m): m \in Z_2\} =
        \{\pi(m): m \in Z_1\cap Z_2\}
    $.
    Therefore, $|\pi(Z_1) \cap  \pi(Z_2)| = |Z_1\cap Z_2|$, which concludes the proof of the first inequality.

    To prove the second inequality, we note that if
    $B_1 \subseteq \bar B_1$ and $B_2 \subseteq \bar B_2$, then
    $ \bar B_1\cap \bar B_2 \subseteq (B_1\cap B_2) \cup (\bar B_1 \setminus B_1) \cup (\bar B_2 \setminus B_2)$,
    and, consequently,
    $$
        |\bar B_1\cap \bar B_2| \le |B_1\cap B_2| + |\bar B_1 \setminus B_1| + |\bar B_2 \setminus B_2|= |B_1\cap B_2| + (|\bar B_1 |-| B_1|) + (|\bar B_2 |-| B_2|).$$
    Hence, due to the assumptions, we obtain
    \[
        \vert A_1 \cap A_2\vert \le \vert \pi(Z_1) \cap  \pi(Z_2) \vert
        +(\vert A_1 |-|  \pi(Z_1) \vert) + (\vert A_2 |-|  \pi(Z_2) \vert).
    \]
    Since $|\pi(Z_1) \cap  \pi(Z_2)|=|Z_1\cap Z_2|$ and  $|\pi(Z_i)| = |Z_i|$, we obtain the second inequality.
\end{proof}
\begin{proof}[Proof of Theorem~\ref{T:1.3}]
    By Lemma~\ref{lemma 1.2},
    \begin{equation*}
        |\pi(i) -\pi(j)|
        =
        l^{(\pi)}_i+ l^{(\pi)}_j + 1 - \left |\left (\pi(i)+ C(l_i^{(\pi)})\right ) \cap \left (\pi(j)+ C(l_j^{(\pi)})\right )\right|.
    \end{equation*}

    By Lemma~\ref{lemma 2.1},
    \[
        \begin{aligned}
            |D_i \cap D_j|
             & \le \left |\left (\pi(i)+ C(l_i^{(\pi)})\right ) \cap \left (\pi(j)+ C(l_j^{(\pi)})\right)\right | \\
             & \le |D_i \cap D_j| + 2l_i^{(\pi)} - |D_i| + 2l_j^{(\pi)} - |D_j|+2.
        \end{aligned}
    \]

    Using the above relations, we obtain upper ($U$) and lower ($L$) bounds
    \begin{multline*}
        |D_i| +  |D_j|   - |D_i \cap D_j| - l_i^{(\pi)} - l_j^{(\pi)}-1 \le |\pi(i)- \pi(j)| \le l^{(\pi)}_i+ l^{(\pi)}_j  - |D_i \cap D_j| + 1.
    \end{multline*}
    Simple algebra leads to the following
    \begin{equation*}
        -\frac{U - L}{2}  \le |\pi(i)- \pi(j)|- \frac{U+L}{2}\le \frac{U - L}{2}.
    \end{equation*}
    Moreover, ${U+ L} = |D_i\triangle D_j| = 2(\mathbf A_{\boldsymbol \Sigma})_{ij}$ and $U - L = 2l^{(\pi)}_i+2l^{(\pi)}_j+2-|D_i| - |D_j|$.
    This ends the proof.
\end{proof}

\begin{proof}[Proof of Proposition~\ref{Prop:2}]
    We apply the inequalities \eqref{eq: approx2} to the formula \eqref{matrix_dist}. Furthermore, since $\mathbf {G}_\pi$ and $\mathbf A_{\boldsymbol{\Sigma}}$ are distance matrices, hence $[\mathbf {G}_\pi]_{ii}=[\mathbf A_{\boldsymbol{\Sigma}}]_{ii}=0$  and
    \begin{align*}
        \delta_m(\mathbf {G}_\pi,\mathbf A_{\boldsymbol{\Sigma}})
         & \le \frac{1}{dm}
        \sum_{\substack{i\ne j                                                                                                                                  \\
        [\mathbf G_\pi]_{ij} \le m}}
        \left(\frac{2l^{(\pi)}_i+1-|D_i|}{2}+\frac{2l^{(\pi)}_j+1-|D_j|}{2}\right)                                                                              \\
         & =\frac{1}{dm}\sum_{i=1}^d\left(2l^{(\pi)}_i+1-|D_i|\right) \Big|\left(\{1,\dots,d\}\cap \left(\pi(i) + C(m)\right)\right) \setminus \{\pi(i)\} \Big| \\
         & =\frac{2}{d}\sum_{i=1}^d(2l^{(\pi)}_i+1-|D_i|) - \frac{1}d
        \sum_{\pi(i)=1}^m \frac{m+1-\pi(i)}{m}(2l^{(\pi)}_{i}+1-|D_{i}|)                                                                                        \\
         & \hphantom{= } -
        \frac{1}d
        \sum_{\pi(i)=d-m+1}^{d} \frac{m - d +\pi(i)}m(2l^{(\pi)}_{i}+1-|D_{i}|).
    \end{align*}
    In the last equation, we used a fact that for $\pi(i) \in \{1,\dots,m\} \cup \{d-m+1,\dots,d\}$ the number of neighbors linearly decreases as $\pi(i)$'s approach the endpoints 1 and $d$. When $m=d$ then the negative terms sum up to $-\frac{d+1}{d^2}\sum_{i=1}^d(2l^{(\pi)}_i+1-|D_i|)$ which together with the positive term give the formula (\ref{eq:uperr bound 1}) after simple transformations.
\end{proof}

For a symmetric $\mathbf A =\left[\mathbf A_{i, j}\right]$, $ l_i(\mathbf A) $ is the half-width of $ D_i $ under the identity permutation.
To prove Theorem~\ref{t:bandopt},
we begin with a lemma describing how the half-width changes when two matrices differ by a single nonzero entry and its symmetric counterpart.
The lemma also characterizes the impact of zeroing out the elements farthest from the diagonal in the original matrix. A rather straightforward proof is omitted.
\begin{lemma}
    Let $\lambda$ be the half-bandwidth of a symmetric $d\times d$ matrix $ \mathbf A $. For $i\in \{1,\dots, d-1\}$ and $k\in\{1,\dots, d-i\}$, let $\boldsymbol A_{i,i+k} $ and  $\boldsymbol A_{i+k,i} $ be nonzero and $ \mathbf B $ be obtained from $ \mathbf A $ by setting these entries to zero. Then
    \begin{equation}
        \|\boldsymbol{l}(\mathbf A)\|_1 - \|\boldsymbol{l}(\mathbf B)\|_1 =  l_i(\mathbf A) - l_i(\mathbf B) + l_{i+k}(\mathbf A) - l_{i+k}(\mathbf B)\ge 0.
        \label{eq:half-widths}
    \end{equation}
    If  $k=\lambda$, then
    \begin{itemize}
        \item[\it i)]  $l_i(\mathbf A) - l_i(\mathbf B) \geq 1$, whenever $\boldsymbol A_{i,i-\lambda} = 0$ or $i\le \lambda$, with the equality holding if, additionally,  $\boldsymbol A_{i,i+\lambda-1} = 1$ or $\boldsymbol A_{i,i-\lambda+1}  = 1$,
        \item[\it ii)] $l_{i+\lambda}(\mathbf A) - l_{i+\lambda}(\mathbf B) \geq 1$, whenever $\boldsymbol A_{i+\lambda, i+2\lambda}  = 0$ or $i+2\lambda>d$, with the equality holding if, additionally,
              $\boldsymbol  A_{i+\lambda, i+1}=1$ or $\boldsymbol A_{i+\lambda,i+2\lambda-1} = 1$.
    \end{itemize}
    \label{l:zeroing one el.}
\end{lemma}

For a symmetric $0$-$1$ matrix $ \mathbf A $ with the half-bandwidth $ \lambda \geq 1 $, let $j = \min \{ i \mid l_i(\mathbf A) = \lambda \}$ and $m(\mathbf A)$ be the largest positive integer $m$ for which
$
    \boldsymbol A_{j + (t-1)\lambda, j + t\lambda}=1,\,\, t=1,\dots, m.
$
It is clear that either $ \boldsymbol A_{j + m\lambda, j + (m+1)\lambda} = 0 $ or $ j + (m(\mathbf A)+1)\lambda > d $.

The sequence of pairs of indices
$$
    S(\mathbf A)=\left(j,j\right)+\lambda \cdot ((0,1),(1,2),\dots ,(m(\mathbf A)-1,m(\mathbf A)) )
$$
describes the first longest contiguous sequence of $1$ along the band, with each element corresponding to an entry in $ \mathbf A $ spaced exactly $ \lambda $ positions apart in row (and column) indices, terminating at either a zero or the matrix boundary.

\begin{corollary}
    Let $ \mathbf A $ have the half-bandwidth $\lambda\geq 1$. Consider the matrix $ \mathbf B $ obtained from $ \mathbf A $ by zeroing out all $1$'s on the boundary given through $ S(\mathbf A) $ together with their symmetric counterparts. Then
    \[
        \|\boldsymbol{l}(\mathbf A)\|_1 - \|\boldsymbol{l}(\mathbf B)\|_1 \geq  m(\mathbf A) + 1
    \]
    and the equality holds if for each $(i,j) \in S(\mathbf A)$, $\boldsymbol A_{i, j-1} = \boldsymbol A_{i+1,j} = 1$.
    \label{cor:BoundaryChain}
\end{corollary}
\begin{proof}
    The proof reduces to iterative application of Lemma~\ref{l:zeroing one el.} for $k=\lambda$. Let  the initial matrix be $ \mathbf A^{(0)} = \mathbf A $. At each step $ t $, $ \mathbf A^{(t)} $ is obtained from $ \mathbf A^{(t-1)} $ by zeroing the entry indexed by the last element in $ S(\mathbf A^{(t-1)}) $ (and its symmetric counterpart).
    Clearly, $ S(\mathbf A^{(t)}) $ consists of the first $ m - t $ elements of $ S(\mathbf A^{(0)}) $, thus, $ \mathbf A^{(m)} = \mathbf B $.

    Due to the construction of the sets $ \mathbf  S(\boldsymbol A^{(t)}) $, it is straightforward that for all $ t < m $, the weaker assumption in {\it ii)} of Lemma~\ref{l:zeroing one el.} is satisfied for the element indexed by the last entry in $ \mathbf  S(\boldsymbol A^{(t)}) $.
    Therefore, using \eqref{eq:half-widths}, we reduce the half-width at each step by at least 1.
    In the case when $ t = m $, the weaker assumptions of {\it i)} and {\it ii)} hold, resulting in a reduction by at least 2.
    When for each $(i,j) \in S(\mathbf A)$, $\boldsymbol A_{i, j-1} = \boldsymbol A_{i+1,j} = 1$, then the stronger assumptions of {\it i)} and {\it ii)} hold at each step, and we obtain the equality.
\end{proof}

\begin{proof}
    [Proof of Theorem~\ref{t:bandopt}]
    Corollary~\ref{cor:BoundaryChain}, applied iteratively, allows us to establish a lower bound for the half-width of any symmetric matrix $\boldsymbol  \Sigma $. Let $ \boldsymbol \Sigma^{(0)} = \boldsymbol \Sigma $.
    At step $ t $, $\boldsymbol \Sigma^{(t)}$ corresponds to $\mathbf B$ in the corollary in which we take $\mathbf A = \boldsymbol \Sigma^{(t-1)}$.
    It is clear that the procedure can be carried out until there are no nonzero entries outside the diagonal and this terminating step is denoted by $\tilde t$, for which $\|\boldsymbol{l}(\boldsymbol \Sigma^{(\tilde{t})})\|_1 = 0$. Thus
    \begin{equation}
        \|\boldsymbol{l}(\boldsymbol \Sigma^{(0)})\|_1
        \geq \tilde{t} + \sum_{i=0}^{\tilde{t}-1} m(\boldsymbol \Sigma^{(i)}) = \tilde{t} + \left |\{ (i,j):\boldsymbol \Sigma_{ij} \neq 0, i<j\}\right |,
        \label{ineq n. iter}
    \end{equation}
    since the iterative procedure eliminates all off-diagonal elements of $\boldsymbol \Sigma^{(0)} $. Furthermore, if $\boldsymbol \Sigma^{(0)} $ is a full band matrix, then
    \begin{equation}
        \|\boldsymbol{l}(\boldsymbol \Sigma^{(0)})\|_1 = \tilde{t} + \left |\{ (i,j):\boldsymbol \Sigma_{ij} \neq 0, i<j\}\right |
        \label{eq n. iter}
    \end{equation}
    since the additional assumptions in {\it i)} and {\it ii)} of Lemma~\ref{l:zeroing one el.} are always fulfilled.

    From \eqref{ineq n. iter}, the lower bound on the half-width is determined by two factors: (1)  the number of nonzero off-diagonal entries and (2) the number of iterations needed to diagonalize the matrix.
    The key to optimality of a permutation is the iteration count, as row and column permutations preserve the number of nonzero elements.

    The number of iterations $\tilde t$ is minimized when the sequences $ S(\boldsymbol \Sigma^{(i)}) $ are as long as possible. Their length depends on two parameters: the half-bandwidth $ \lambda_i $ and the number of positions on the subdiagonal at order $ \lambda_i $, which equals $d - \lambda_i$.  A larger $\lambda_i$ increases spacing between consecutive elements  in $ S(\boldsymbol \Sigma^{(i)}) $ and reduces available positions, making $m(\boldsymbol{\Sigma}^{(i)})$ a decreasing function of $\lambda_i$.
    From this reasoning, combined with \eqref{eq n. iter}, we conclude that among all matrices with the number of nonzero off-diagonal elements equal to the number in a full band matrix, the full band matrix configuration requires the minimal number of iterations.
    This proves Theorem~\ref{t:bandopt}.
\end{proof}

\begin{proof}[Proof of Theorem~\ref{T: per full band matr}]
    The part {\it i)} is a rather standard observation due to the connection between $\boldsymbol \Sigma(t)$ and the matrix $\boldsymbol \Sigma^t$ and noticing that for a matrix $\mathbf A$ and a permutation matrix $\tilde{\boldsymbol \pi}$,
    $
        \left(\tilde{\boldsymbol{\pi}}^{\top} \mathbf A\tilde{\boldsymbol{\pi}}\right)^t=\tilde{\boldsymbol{\pi}}^{\top} \mathbf A^t\tilde{\boldsymbol{\pi}}.
    $

    In order to prove {\it  ii)} and {\it
            iii)}, it is enough to consider properties of  a $ d \times d $ full band matrix with half-bandwidth $ \lambda \geq 1 $.
    Under {\it ii)}, only two rows in the original (unpermuted) band matrix have $\lambda + 1$ nonzero entries—the first and the last—yielding two choices for the first row. The remaining rows must then follow uniquely: in the original order if the first row is selected, or in reverse if the last row is chosen first.

    To prove {\it iii)}, observe that, by the same argument as before, there are only two ways to choose the first and last $d - \lambda - 1$ rows. The remaining $2\lambda + 2 - d$ central rows are fully nonzero, so their permutations do not affect the matrix's sparsity pattern. The final relation follows from the fact that the full central rows that can be arbitrarily permuted for $\boldsymbol{\Sigma}(t)$ constitute a subset of those that correspond to $\boldsymbol \Sigma(t+1)$.

\end{proof}

\section{Classical multidimensional scaling}
\label{CMDS_supp}

Here we present an overview of the fundamental concepts behind the classical multidimensional scaling (CMDS) method.
The method is also  referred to as Principal Coordinates Analysis, Torgerson  (or Torgerson–Gower) Scaling, \cite{Torgerson52, GOWER66}.
CMDS aims to find a configuration of $d$ points in a $ p $-dimensional Euclidean space represented as rows in $ \mathbf{X} \in \mathbb{R}^{d \times p} $, $d>p$,  such that distances between them approximate the original distances given in the $ d \times d $ distance matrix $ \mathbf{G} $, see Algorithm~\ref{alg:MDS_supp}.
\begin{algorithm}
    \caption{Classical multidimensional scaling}
    \label{alg:MDS_supp}
    \begin{algorithmic}\vspace{2mm}
        \STATE{{\sc Input:} $\mathbf G = [G_{ij}]$ - $d\times d$ symmetric matrix of  distances;
        $p<d$ -- reduced dimension }\vspace{1mm}
        \STATE{$ \mathbf{K} = \begin{bmatrix}G_{ij}^2\end{bmatrix} $};  \COMMENT{{\tt \small Compute the squared distance matrix}}\vspace{2mm}
        \STATE{$ \mathbf{H} = \mathbf{I} - \frac{1}{d} \mathbf{1}\mathbf{1}^\top $}; \COMMENT{{\tt \small 1 - the column vector of 1's }}\vspace{2mm}
        \STATE{$
                \mathbf{B} = - \mathbf{H} \mathbf K \mathbf{H}/2.
            $}; \COMMENT{{\tt \small Double-centering}}\vspace{2mm}
        \STATE{$
            (\boldsymbol \Lambda, \mathbf{V}) = {\rm eigen}(\mathbf B)
        $}; \COMMENT{{\tt \footnotesize Decreasing eigenvalues, and the corresponding eigenvectors}}\vspace{2mm}
        \STATE{$
                \boldsymbol \Lambda_p= \boldsymbol \Lambda[1:p]$};  \COMMENT{{\tt \small The p-largest eigenvalues, assumed to be positive}}\vspace{2mm}
        \STATE{$
                \mathbf V_p= \mathbf V[,1:p]$};  \COMMENT{{\tt \small The corresponding eigenvectors}}\vspace{3mm}
        \STATE{{\sc Return:} $\mathbf X_p=\mathbf V_p\, {\rm diag}(\boldsymbol \Lambda_p)^{1/2}$ - The $d$ rows are $p$-dimensional configuration vectors with distances approximating the ones in $\mathbf G$.  }
    \end{algorithmic}
\end{algorithm}

The rationale behind the algorithm can be easily explained in the special case when the initial matrix $\mathbf G$ is made of the Euclidean distances, i.e., the distances are between vectors that are rows in a $d\times d$ matrix $\mathbf Y$, i.e., $\mathbf  G_{ij}=\|\mathbf Y_{i\cdot}-\mathbf Y_{j\cdot}\|$  (additionally, let  $\mathbf K = [\mathbf G_{ij}^2]$).
Clearly, for $\boldsymbol \Sigma= \mathbf Y \mathbf Y^\top$
and $d\times 1$ matrix $\mathbf 1$ of ones
\begin{equation}
    \label{eq:G}
    \mathbf K=\mathbf 1\,{\rm diag}(\boldsymbol \Sigma)^\top +{\rm diag}(\boldsymbol \Sigma)\, \mathbf 1^\top -2 \boldsymbol \Sigma.
\end{equation}
Here, ${\rm diag}(\cdot)$ extracts the diagonal vector from a matrix argument, or creates a diagonal matrix from a vector argument.
If each of the columns is recentered by a constant, then $\mathbf G$ does not change.
Consequently, we may assume that the columns of $\mathbf Y$ are centered at zero, i.e., $\boldsymbol \Sigma \mathbf 1=\mathbf 1^\top \boldsymbol \Sigma=0$.
We observe also that $\mathbf H \mathbf 1=\mathbf 1^\top \mathbf H=0$ so that, with the notation of Algorithm~\ref{alg:MDS_supp}, we have
\begin{align*}
    \mathbf H \mathbf K \mathbf H & = \mathbf H\, \mathbf 1\,{\rm diag}(\boldsymbol \Sigma)^\top \mathbf H
    +
    \mathbf H\,{\rm diag}(\boldsymbol \Sigma)\, \mathbf 1^\top \mathbf H
    -2 \mathbf H\boldsymbol \Sigma \mathbf H                                                               \\
                                  & =
    -2 \boldsymbol \Sigma +\frac{2}d \left(\mathbf 1 \mathbf 1 ^\top\boldsymbol \Sigma +   \boldsymbol \Sigma  \mathbf 1\mathbf 1 ^\top \right) - \frac 2{d^2} \mathbf 1 \mathbf 1^T \boldsymbol \Sigma \, \mathbf 1 \mathbf 1^\top =-2 \boldsymbol \Sigma.
\end{align*}
From this, we observe that $\mathbf B$ in the algorithm is equal to positive semidefinite $\boldsymbol \Sigma$ and since $\boldsymbol \Sigma=\mathbf V \boldsymbol \Lambda \mathbf V^\top$ thus $\boldsymbol \Sigma - \mathbf X_p\mathbf X_p^\top$ is also positive semidefinite having the sum of eigenvalues (trace norm) equal to $\lambda_{p+1}+\dots +\lambda_d$, where $\boldsymbol \Lambda = \operatorname{diag}(\lambda_1,\dots, \lambda_d)$.
Consequently, replacing in \eqref{eq:G}, $\boldsymbol \Sigma$ by $\mathbf X_p\mathbf X_p^\top$ leads to a good approximation $\mathbf G_p=\|\mathbf X_{p;i\cdot}-\mathbf X_{p;j\cdot}\|$ of $\mathbf G$ with the error controlled by $\lambda_{p+1}+\dots +\lambda_d$.
When $\mathbf G$ is not generated by Euclidean distances, the error analysis is no longer straightforward.
For a more detailed and comprehensive description of the method, readers are referred to \cite{BorgGroenen05}.

The most computationally intensive part of CMDS is the eigendecomposition of the $ d \times d $ matrix $ \mathbf{B} $, which has a complexity of $ O(d^3) $ (for details see e.g. \cite{trefethen_numerical_1997}).
Additionally, storing the $ \mathbf K $ and $ \mathbf{B} $ matrices requires $ O(d^2) $ memory.
In Section~\ref{sec: loc_config} of the paper,  we show that CMDS is used in our method to locally order the sets of nearest neighbors $D_i$. Since we consider sparse matrices, one can expect that the cardinalities $|D_i|$ are relatively small, and CMDS can therefore find a solution in a relatively short time. However, in cases where $|D_i|$ are large, various algorithms exist that are capable of efficiently handling even large distance matrices. For an accessible overview of such methods, see \cite{Delicado2024}.

\section{Retrieving ordering from distance matrices}
\label{DistMat_supp}
Let us consider the distance matrix $\mathbf G_{\mathbf x}$ for points on the real line as defined at the beginning of Section~\ref{sec: dist_matrix} of the paper.
The algorithm that reconstructs $\mathbf x$ using only the information in $\mathcal I(\mathbf G_{\mathbf x})$ serves as the proof of the following theorem.

\begin{theorem}
    \label{th:config}
    Let $\mathbf{x} \in \mathbb{R}^d$ be such that $x_i \ne x_j$ for all $i \ne j$, and define the matrix $\mathbf{G}_{\mathbf{x}} = \left[\,|x_i - x_j|\,\right]_{i,j = 1}^d$.
    Then there exists $\mathbf y \in \mathbb R^d$ that depends on $\mathbf G_{\mathbf x}$ only through $\mathcal I(\mathbf G_{\mathbf x})$ such that $\mathbf G_{\mathbf y}=\mathbf G_{\mathbf x}$.  Moreover, there exists $s\in \mathbb R$ such that either $\mathbf x=s+\mathbf y$ or $\mathbf x=s-\mathbf y$.
\end{theorem}
\begin{proof}

    Let permutation $\beta$ sort $x_i$, i.e., $x_{\beta(1)}<\dots < x_{\beta(d)}$.
    Given $\mathcal I(\mathbf G_{\mathbf x})$, we choose one of the two endpoint indices, e.g., $\alpha(1)=j_e^1$, and set $y_1=0$.
    Two cases arise: either $\alpha(1)=\beta(1)$ or $\alpha(1)=\beta(d)$.
    We begin with the former.

    Furthermore, there is only one $i\in \{1,\dots, d\}\setminus \{\alpha(1)\}$ such that $\alpha(1)\in \{j_i^1,j_i^2\}$.
    We define $\alpha(2)=i$ and $y_2=\rho$, where $\rho$ is the one from $\{\rho_i^1,\rho_i^2\}$ that corresponds to $\alpha(1)$.
    We observe that $\alpha(2)=\beta(2)$ and
    $|y_2-y_1|=y_2-y_1=x_{\alpha(2)}-x_{\alpha(1)}=|x_{\alpha(2)}-x_{\alpha(1)}|.
    $
    Then we define recursively $\alpha(k)$,  and $y_k$, for $k=3,\dots, d-1$.\\
    Given $\alpha(1),\dots, \alpha(k-1)$ and $y_1<\dots <y_{k-1}$, such that
    \begin{equation}
        \label{eq:inc}
        \begin{split}
            \alpha(l)                      & =\beta(l),\,\, l=1,\dots , k-1,   \\
            x_{\alpha(l+1)}-x_{\alpha(l)}= & y_{l+1}-y_l,\,\, l=1,\dots , k-2,
        \end{split}
    \end{equation}
    we define $\alpha(k)$ as the only $i \notin \{\alpha(1),\dots, \alpha(k-1)\} $ such that $\alpha(k-1)\in  \{j_i^1,j_i^2\}$. Furthermore, we define $y_k=y_{k-1}+\rho$, where $\rho$ is the one from $\{\rho_i^1,\rho_i^2\}$ that corresponds to $\alpha(k-1)$. Clearly, $x_{\alpha(k)}$  corresponds to the in-between point associated with $\{(j_i^1, \rho_i^1),(j_i^2, \rho_i^2)\}$ and thus $\alpha(k)=\beta(k)$.
    This ensures that \eqref{eq:inc} is also satisfied for $l=k-1$.
    Moreover, $\alpha(d)$ is the only remaining index to the sequence $\alpha$ and $y_d = y_{d-1} + \rho$, where $\rho$ is the one from $\{\rho_{\alpha(d-1)}^1,\rho_{\alpha(d-1)}^2\}$ that corresponds to $\alpha(d)$. Clearly, \eqref{eq:inc} remains satisfied for
    $k=d$ as well.

    Redefine $\mathbf y$ as $y_{\alpha^{-1}(k)}$, $k=1,\dots, d$. By \eqref{eq:inc},
    $
        x_{\alpha(l+1)}-x_{\alpha(l)}=y_{\alpha(l+1)}-y_{\alpha(l)},
    $
    and thus, for $s=x_{\alpha(1)}$,
    $
        \mathbf x = s+\mathbf y.
    $
    Consequently, $\mathbf G_x = \mathbf G_y$.

    For the case when $\alpha(1)=\beta(d)$, consider $\tilde{\mathbf x}=-\mathbf x$ for which $\alpha(1)=\tilde \beta(1)$ and $\mathbf G_{\mathbf x}=\mathbf G_{\tilde{\mathbf x}}$. By the previous argument we construct ${\mathbf y}$ such that
    $
        \tilde{\mathbf x} =-\mathbf x = \tilde s+\mathbf y$, which means that $\mathbf x=-\tilde s-\mathbf y$.
\end{proof}

\section{Examples of the R Package \texttt{DyadiCarma}}

The \texttt{DyadiCarma} package is available on the Comprehensive R Archive Network (CRAN). It provides support for dyadic matrix objects, including arithmetic operations, conversion between standard R matrices and dyadic matrix objects, and the dyadic factorization algorithm introduced in the main text. The package relies on \texttt{Rcpp} and \texttt{RcppArmadillo} for computational efficiency. In this section, we present two illustrative examples demonstrating the use of \texttt{DyadiCarma}.

\begin{example}[Arithmetic operations of dyadic matrices]
    \begin{lstlisting}[language=R, breaklines=true, basicstyle=\small\ttfamily]

    library(DyadiCarma)
    N <- 4; k <- 3

    V  <- construct(N, k, type = "vert")
    H  <- construct(N, k, type = "horiz")
    S  <- construct(N, k, type = "symm")
    AS <- construct(N, k, type = "asymm")

    # Linear combination of dyadic matrices
    linearComb <- -S + 3 * H - 6 * AS + V

    # Matrix conversion for identity check
    identical(
    as.matrix(linearComb),
    -as.matrix(S) + 3 * as.matrix(H) - 6 * as.matrix(AS) + as.matrix(V)
    )
    ## [1] TRUE

    # Matrix type
    linearComb@type
    ## [1] "asymm"
\end{lstlisting}
\end{example}

\begin{example}[Inversion of a dyadic matrix with the dyadic algorithm]
    \begin{lstlisting}[language=R, breaklines=true, basicstyle=\small\ttfamily]

    library(DyadiCarma)
    N <- 4; k <- 3

    # Construct a 45x45 vertically dyadic matrix with random entries
    V <- construct(N, k, type = "vert", distr = "unif")

    # Construct a 45x45 symmetric dyadic matrix by V^T V
    S <- t(V) %*% V
    S@type <- "symm"
    S@aentries <- list() # Convert S from "asymm" to "symm"

    # Find the vertically dyadic matrix P such that t(P) %*% S %*% P = I
    # using the dyadic factorization algorithm.
    # The three "max" results below may differ slightly due to randomness.

    P <- dyadFac(S)
    I1 <- as.matrix(t(P) %*% S %*% P)
    I  <- diag(dim(I1)[1])
    max(abs(I1 - I))
    ## [1] 9.624634e-12

    # Compute the inverse of S using the dyadic algorithm
    matS <- as.matrix(S)
    iS <- dyadFac(S, inv = TRUE)
    I2 <- iS %*% matS
    max(abs(I2 - I))
    ## [1] 9.924609e-11

    # Compare with the result using the built-in solve() function
    iS_solve <- solve(matS)
    I3 <- iS_solve %*% matS
    max(abs(I3 - I))
    ## [1] 2.257426e-08
\end{lstlisting}

\end{example}

\end{document}